\documentclass{article}

\usepackage{geometry}                % See geometry.pdf to learn the layout options. There are lots.
\usepackage{graphicx}

\usepackage{epstopdf}

\usepackage{enumerate}

\usepackage{amsmath}
\usepackage{amssymb}
\usepackage{amsthm}
\usepackage{url}
\usepackage{color}

\usepackage{tikz-cd}

\usepackage{float}
\floatstyle{boxed} 
\restylefloat{figure}
\usepackage[section]{placeins}

\usepackage{chngcntr}
\counterwithin{figure}{section}

\usepackage{pifont}% http://ctan.org/pkg/pifont
\usepackage{stmaryrd}
\usepackage{dsfont}
\usepackage{textcomp}
\usepackage{verbatim}
\usepackage{accents}
\usepackage{wasysym}
\usepackage{csquotes}
\usepackage{booktabs}  % For better-looking tables

\newtheorem{thm}{Theorem}
\newtheorem*{thm*}{Theorem}
\newtheorem{prop}[thm]{Proposition}
\newtheorem{lemma}[thm]{Lemma}

\newtheorem*{mainthm*}{Main Theorem}
\newtheorem*{mainlemma*}{Main Lemma}
\newtheorem{cor}[thm]{Corollary}

\newtheorem*{claim*}{Claim}

\newtheorem*{conj*}{Conjecture}

\theoremstyle{definition}

\newtheorem{dfn}[thm]{Definition}
\newtheorem*{dfn*}{Definition}

\newtheorem*{fix*}{Fixation}

\newtheorem*{tgtfeat*}{Target features}

\newtheorem*{ass*}{Assumption}

\newtheorem*{que*}{Question}

\newtheorem*{chal*}{Challenge}

\newtheorem*{aim*}{Aim}

\newtheorem*{diag*}{Diagnosis}
\newtheorem*{mque*}{Main Question}

\newtheorem*{rem*}{Remark}

\newtheorem*{appr*}{Approach}
\newtheorem*{Copernican*}{Copernican Principle}
\newtheorem*{Correctness*}{Correctness Condition}

\newcommand{\Lang}{\mathcal{L}}
\newcommand{\Atom}{\mathrm{Atom}}

\newcommand{\Term}{\mathrm{Term}}

\newcommand{\ev}{\mathrm{ev}}
\newcommand{\Var}{\mathrm{Var}}

\newcommand{\dom}{\mathrm{dom}}

\newcommand{\udot}[1]{\text{\d{${#1}$}}}

\newcommand{\gquote}[1]{\ulcorner #1 \urcorner}

\newcommand{\lrangle}[1]{\langle #1 \rangle}

\newcommand{\ar}{\mathrm{ar}}
\newcommand{\lng}{\mathrm{length}}

\newcommand{\sbt}{\mathrm{sbt}}
\newcommand{\gq}{\mathrm{gq}}

\newcommand{\df}{\mathrm{df}}

\newcommand{\dt}{\hspace{2pt}}

\newcommand{\PA}{\mathsf{PA}}

\newcommand{\ZF}{\mathsf{ZF}}
\newcommand{\ZFC}{\mathsf{ZFC}}

\newcommand{\UCT}{\mathsf{UCT}}
\newcommand{\FS}{\mathsf{FS}}

\newcommand{\NEC}{\mathsf{NEC}}
\newcommand{\CONEC}{\mathsf{CONEC}}

\usepackage[square,numbers]{natbib}

\newcommand{\statsub}{\subseteq_{\Lang_\mathrm{StatPred}}}
\newcommand{\dynin}{\mathrel{\in^*}}

\newcommand{\tildin}{\mathrel{\tilde{\in}}}

\title{Development Processes}
\author{
Paul Gorbow \\
}
\date{}                                           % 

\begin{document}

\maketitle

\begin{abstract}
Throughout mathematics there are constructions where an object is obtained as a limit of an infinite sequence. Typically, the objects in the sequence improve as the sequence progresses, and the ideal is reached at the limit. I introduce a view that understands this as a development process by which a dynamic mathematical object develops teleologically. In particular, this paper elaborates and clarifies the intuition that such constructions operate on a single dynamic object that maintains its identity throughout the process, and that each step consists in a transformation of this dynamic object, rather than in a genesis of an entirely new static object. This view is supported by a general philosophical discussion, and by a formal modal first-order framework of development processes. In order to exhibit the ubiquity of such processes in mathematics, and showcase the advantages of this view, the framework is applied to wide range of examples: The set of real numbers, forcing extensions of models of set theory, non-standard numbers of arithmetic, the reflection theorem schema of set theory, and the revision semantics of truth. Thus, the view proposed promises to yield a unified dynamic ontology for infinitary mathematics. 
\end{abstract}

\section{Introduction}\label{Sec:Intro}

Consider a process in which some entity develops to satisfy certain significant properties. I call this type of process a {\it development process}, and the entity in question a {\it  dynamic entity}. Note that development processes have a dynamic and teleological nature. For example, the Leibniz series 
\[
\frac{4}{1} - \frac{4}{3} + \frac{4}{5} - \frac{4}{7} + \cdots \approx \pi
\]
may be viewed as a development process in which a dynamic entity, namely a rational approximation, develops to gain certain significant properties, namely successively smaller approximation error to $\pi$. The process has the structure of a sequence of states, in which the rational approximation, say $\alpha$, takes on different values in different states, $s$: 

\[
\renewcommand{\arraystretch}{2} % Adjusts the space between rows
\begin{array}{c|c|c|c|c|c}
\hline
s & 0 & 1 & 2 & 3 & \cdots \\ \hline
\alpha  &  4  &  4 - \frac{4}{3}  &  4 - \frac{4}{3} + \frac{4}{5}  &  4 - \frac{4}{3} + \frac{4}{5} - \frac{4}{7}    &  \cdots  \\ \hline
\textrm{error bound} & \frac{4}{3}  &  \frac{4}{5}  &  \frac{4}{7}  &  \frac{4}{9}  &  \cdots \\ \hline
\end{array}
\]

This paper develops a philosophical account and a formal framework for development processes in mathematics, and carefully applies it to a diverse array of mathematical phenomena. First, I characterize development processes and dynamic entities philosophically in \S \ref{Sec:Char_Dev_Proc}. Then I proceed, in \S \ref{Sec:Prel}, to formalize development processes with the tools of first-order logic with modal operators. As shown in \S \ref{Sec:Background}, this formal framework is a generalization of the modal operator approach to potentialism from \cite{Lin13, LS19}. Its relationships to Carnap's individual concepts \cite{Car56}, to divergent potentialism and free choice sequences \cite{BLS22, Bra23}, and to arbitrary objects \cite{Fin85a, Fin85b, HS19, VY24}, are also explained. \S \ref{Sec:Tel_Dev} provides mathematical results on the teleological aspect of development processes in general. The ensuing sections are devoted to the following applications:
\begin{description}
\item[\textnormal{ \S \ref{Sec:Reals}}] The set of real numbers may be represented by a development process on a set of dynamic rational numbers with the appropriate Cauchy convergence property. %(\S \ref{Sec:Reals})
\item[\textnormal{ \S \ref{Sec:Forcing}}] A forcing extension in set theory may be viewed as development processes in which the $\in$-relation develops to satisfy certain properties. %(\S \ref{Sec:Forcing})
\item[\textnormal{ \S \ref{Sec:Saturation}}] A non-standard number, appearing in an ultrapower of the standard model $\mathbb{N}$ of the natural numbers, may be represented by a development process on a dynamic natural number. This generalizes to realizers of model-theoretic types quite generally. %(\S \ref{Sec:Saturation})
\item[\textnormal{ \S \ref{Sec:Reflection}}] The Reflection theorem schema of set theory may be given a very simple formulation as a theorem schema of a development process. %(\S \ref{Sec:Reflection})
\item[\textnormal{ \S \ref{Sec:Revision_semantics}}] The class of true sentences of arithmetic, in a language with an untyped truth predicate, may be represented by the revision semantics, construed as a development process on a dynamic set or class of sentences. %(\S \ref{Sec:Revision_semantics})
\end{description}
Finally, the concluding \S \ref{Sec:Conc} summarizes the applications, and drives the thesis that the dynamic perspective of development processes has a unifying power for our philosophical understanding of infinitary mathematics.

%This treatment facilitates a balanced potentialist position in the philosophy of mathematics, between {\it  finitism} and {\it  infinitism}: As exemplified by the Leibniz series above, each state of the dynamic object can be obtained by {\it  finitely} many steps, but the process is {\it  infinite} in the sense that there is potential for further development from any state. The dynamic object lives in an ambient structure, such as the rational or natural numbers. In many applications, this ambient structure can in turn be represented by a development process on a dynamic structure that is finite in every state of the process, but converges to an infinite structure.

\section{Development processes and dynamic entities}\label{Sec:Char_Dev_Proc}

This section outlines the general characteristics of the development processes and dynamic entities studied in this paper, as well as some significant distinctions between such processes. %The formalizations introduced in the ensuing sections are intended to precisely articulate and theoretically elucidate such development processes.

There are two common representations of processes whereby a mathematical entity develops through states:
\begin{description}
\item[Function{\normalfont :}] A function (a set of ordered pairs) from the set of states of the process, mapping each state to a value. For the Leibniz series it would be the function $f : \mathbb{N} \rightarrow \mathbb{Q}$, identified with the set $\big\{\langle 0, 4 \rangle, \langle 1, 4 - \frac{4}{3} \rangle, \langle 2,  4 - \frac{4}{3} + \frac{4}{5} \rangle, \cdots \big\}$.
%\item[Intensional function] A formula that defines the above extensional function. For the Leibniz series it could be the formula $f(s) = \Sigma_{i = 0}^s \frac{4}{(-1)^i(2i+1)}$, or the corresponding recursive definition, in the context of the theory of the ring of rational numbers.
\item[Limit{\normalfont :}] A limit value of the process, which abstracts away from the dynamics of the process. For the Leibniz series it would be the irrational number $\pi$.
\end{description}

However, none of these representations captures a single dynamic entity that develops. The function represents the entire process as a collection of distinct values associated with different states, lacking a unified identity across states. The limit represents a static and ideal entity that does not reflect the ongoing dynamics of the process.

To accurately represent a single entity that retains its identity while undergoing change, I introduce the notion of dynamic entities in mathematics. There is an analogy here with temporal objects in metaphysics: just as we may intuitively think of physical objects as persisting through time (endurantism), we can conceive of dynamic mathematical entities as persisting through states of a process. Both realms involve entities that undergo change while maintaining a unified identity. This conception is particularly compelling for processes characterized by development, motivating a framework that can represent a dynamic mathematical entity as it evolves. 

Modal logic, which is well-suited for representing temporal entities in philosophy, is likewise natural for modeling dynamic entities in mathematics. Below, I give a blueprint for a formal structure of development processes with dynamic entities using possible worlds (formal details are provided in \S \ref{Sec:Prel}).

\begin{description}
\item[States{\normalfont :}] A process has a set or class of states, which the process transitions to and from. Each state $s$ is mapped to a possible world $u_s$, which is a first-order structure modeling that state. 
\item[Accessibility relation{\normalfont :}] The states are related by an accessibility relation, which is a partial ordering. This ordering represents the direction of the process. 
\item[Language and semantics{\normalfont :}] $L$ is a fixed first-order language appropriate for describing the static aspects of any state. $L_\mathrm{Dyn}$ is $L$ augmented with so-called {\it  dynamic ideology}, in particular it has a binary relation $\leftleftarrows$ explained below. $L^\Diamond$ and $L_\mathrm{Dyn}^\Diamond$ are $L$ and $L_\mathrm{Dyn}$, respectively, augmented with the modal operators $\Box, \Diamond$. For each state $s$, there is a structure $u_s$, which interprets $L_\mathrm{Dyn}$. For any $\phi \in L^\Diamond$, $s \models \phi$ is defined by the usual possible worlds interpretation of the operators $\Box, \Diamond$.
\item[Static entities{\normalfont :}] The {\it  static entities} are the entities which remain fixed throughout the process. The {\it  static individuals} constitute a sort of individuals.
\item[Dynamic entities{\normalfont :}] The {\it  dynamic entities} are the entities which are developed in the process. The {\it  dynamic individuals} constitute a sort of individuals. $L_\mathrm{Dyn}$ has a binary relation $\leftleftarrows$ from the sort of static individuals to the sort of dynamic individuals. For any state $s$, dynamic individual $\delta$, and static individual $a \in u_s$, $s \models a \leftleftarrows$ expresses that $a$ is {\it  the manifestation of $\delta$ in} $s$. (Or $\delta$ {\it  manifests (as $a$) in} $s$.) This relation is formalized by axioms ensuring that, in each state, each dynamic individual manifests as at most one static individual, and that if it manifests, then it manifests in all accessible states. Analogously, notions of {\it  dynamic classes} and {\it  dynamic structures} will also be formalized.
\item[Type{\normalfont :}] The static individuals will typically share a type $T$.\footnote{Although this paper does not utilize a formal type-theoretic framework, it is useful to employ an informal theoretical notion of type. A type is a type of mathematical object, such as natural number, rational number, or model of the theory of dense linear orders. All objects of a certain type have certain similarities, which make certain relations and operations applicable to all of them. For example, all rational numbers can be considered as the ratio of two integers, and a consequence of this is that a natural operation of multiplication on the rational numbers can be defined in terms of the usual multiplication of integers. } In the case of the Leibniz series they are all of the type rational number. An {\it  entity of type dynamic $T$} (or a {\it  dynamic entity of type $T$}) is a dynamic entity that manifests as an entity of type $T$ in every state where it manifests. For example, the approximation $\alpha$ above is a dynamic rational number, because it manifests as a rational number in each state.
\end{description}

This approach of formalizing development processes in modal logic typically yields systems with tempered logical complexity, because the modal operators have limited expressive power concerning the structure of the accessibility relation on the states. In contrast, the formalization of development processes through functions is typically based on set theory, or at least arithmetic, which is highly complex. For example, if the set of states of a process is the standard model of arithmetic, $\mathbb{N}$, then the typical ambient theory of the functional representation of this process is extremely expressive with regard to its structure. Note that the tempered complexity of the present approach is enabled by the notion of dynamic entity. 

%We can go further and consider what ontological resources are needed in each state. Perhaps the numbers $3$, $4$, and $\frac{4}{3}$ are sufficient in state $0$, while $3$, $4$, $5$, $\frac{4}{3}$, and $\frac{4}{5}$ are sufficient in state $1$. Thus, each state $s$ could instead be mapped to a finite substructure $\mathcal{P}_s$ of $\mathcal{Q}_s$ that contains sufficient ontological resources. It is natural to also consider the sequence of these finite structures as a development process, where a {\it  dynamic structure} is successively extended. There is particular interest in such sequences that {\it  converge to} $\mathbb{Q}$, defined as that for each $q \in \mathbb{Q}$ and each state $s$, there is a state $t$, accessible from $s$, such that $q \in \mathcal{P}_t$.

The nature of dynamic entities possesses three fundamental attributes: type coherence, dynamics, and unity.

\begin{description}
\item[Type coherence{\normalfont :}] A dynamic entity of type $T$ manifests as an entity of type $T$ in every state. 
\item[Dynamics{\normalfont :}] A dynamic entity is inherently dynamic, changing its manifestation from state to state along the directions dictated by the accessibility relation.
\item[Unity{\normalfont :}] Although a dynamic entity of type $T$ may manifest as different entities of type $T$ in different states, it remains the same entity throughout the process.
\end{description}
Thus a dynamic entity evolves through the process while maintaining its unity and type. This allows us to apply the conceptual tools associated with the type $T$ to the dynamic entity. For example, a dynamic rational number is the quotient of two dynamic integers.

In contrast to dynamic entities, traditional representations as functions and limits fail to combine the above attributes: A function mapping states to entities of type $T$ considers these values as plural. The unity provided by the function is of a very different type than $T$. Introducing a limit entity to unify the process results in a purely static representation that abstracts away the very dynamics that characterizes the process. Moreover, it is often of a different type than the dynamic entity (e.g., $\pi$ is not a rational number). The value of introducing dynamic entities (of type $T$) to the ontology of mathematics lies in that they are unitary entities that develop, while retaining their identity through this development, and that the conceptual toolbox endowed by the type $T$ applies to them. Therefore, this approach offers a richer and more accurate representation of mathematical processes, aligning with the familiar endurantist intuition about the evolution of objects.

Here follow some notions applicable to development processes, and how they are formally represented in this paper:

\begin{description}
%\item[Language] The language is used to express the properties of, and relations between, the (dynamic and static) entities.
%\item[States] Each state has a domain of static entities, and determines which formulas in the language are satisfied by which static objects.
%\item[Transition relation] Any state $s$ is related to the states to which the process can transition from $s$ (these are called the accessible states from $s$). Every state accesses itself.
%\item[Values] In any state, each dynamic entity takes on a value. This value is a static entity.
\item[Transformation operations] A {transformation operation} is an operation that can potentially be applied to a dynamic entity to transform it as the process transitions from a state to another state.
\item[Transformations] A {\it  transformation} is a particular application of a transformation operation to a dynamic entity as the process transitions from a state to another state.
\item[Teleological development] The {\it  teleology} of a dynamic entity is a collection of conditions that the dynamic entity develops in order to meet. For any state, and any condition $C$, $C$ is satisfied {\em teleologically} if there is the potential to transition to a state at which $C$ holds for all accessible states.
\item[Non-finality] A development process is {\it  non-final} if there is no state at which the whole teleology is met.
\end{description}

Mathematical objects are often considered immutable, a position emphatically supported by Plato in his theory of forms. There are two prominent understandings of their immutability: one is that they cannot change over time, and another is that they metaphysically could not have been different. But the modality considered in this paper is the modality of change through a conceptual process, which is neither the temporal nor the metaphysical modality. Thus, it is perfectly reasonable for a Platonist to consider modalities in which some mathematical objects are dynamic. Indeed, a mathematical development process, such as an approximation process, may be viewed as a Platonic form containing an internal modality under which the approximation in question changes.

\section{Formalizing development processes }\label{Sec:Prel}

This section formalizes development processes in classical first-order modal logic, with the modalities expressed by the usual operators $\Box, \Diamond$. The primitive connectives, quantifiers and operators are $\neg, \vee, \exists, \Diamond$, and the other ones are defined in terms of these in the standard way.

Throughout the paper, the variable $n$ is assumed to range over the natural numbers. 

Let $K$ be an arbitrary sorted first-order language or sorted first-order language with modal operators, and with a sort called $\mathrm{Stat}$. The sort $\mathrm{Stat}$ is meant to be interpreted by a domain of static individuals. $K$ also denotes the set of $K$-formulas. A {\it  literal} of $K$ is a formula of the form $R(x_1, \cdots, x_n)$ or $\neg R(x_1, \cdots, x_n)$, where $R$ is an $n$-ary relation symbol of $K$ on sort $\mathrm{Stat}^n$. Let $K'$ be a language that is obtained from $K$ by adding some sorts and/or relation symbols. Then $K'$ is called an {\it  expansion} of $K$, and $K$ is called a {\it  reduct} of $K'$. If, moreover, $\mathcal{M}, \mathcal{M}'$ are models interpreting $K, K'$, respectively, such that $\mathcal{M}$ and $\mathcal{M}'$ have the same frames and domains, and interpret $K$ identically, then $\mathcal{M}'$ is called a {\it  ($K'$-)expansion} of $\mathcal{M}$, and $\mathcal{M}$ is called the {\it ($K$-)reduct} of $\mathcal{M}'$. Moreover, the {\it ($K$-)reduct} of $\mathcal{M}'$ is denoted $\mathcal{M}' \restriction K$. If $K$ does not have modal operotors, then $K^\diamond$ denotes the expansion of $K$ by the first-order modal operators $\Box, \Diamond$. 

For a tuple $\vec{t}$, $\lng(\vec{t})$ denotes its length, and $\pi_i(\vec{t})$ denotes the $i$-th component of $\vec{t}$, that is, $\vec{t} = \langle \pi_1(\vec{t}), \cdots, \pi_{\lng(\vec{t})}(\vec{t}) \rangle$. For a relation symbol $R$ or function symbol $f$, $\ar(R)$ and $\ar(f)$ denote their respective arities. Let $\vec{v}$ be a tuple of variables, $R$ be a relation symbol, $f$ be a function symbol, and $S_1, \cdots, S_n$ be sorts. $\vec{v} : (S_1, \cdots, S_n)$ denotes that $\lng(\vec{v}) = n$ and $\pi_1(\vec{v}), \cdots, \pi_n(\vec{v})$ are of the sorts $S_1, \cdots, S_n$, respectively. $R : (S_1, \cdots, S_n)$ denotes that $R$ is an $n$-ary relation symbol that applies to tuples $\vec{t}$, such that $\vec{t} : (S_1, \cdots, S_n)$. $f : (S_1, \cdots, S_n) \rightarrow S$ denotes that $f$ is a function symbol that applies to tuples $\vec{t}$, such that $\vec{t} : (S_1, \cdots, S_n)$ and returns values of sort $S$. If a tuple of variables or a relation symbol is introduced without a sort declaration, all the relevant sorts are assumed to be $\mathrm{Stat}$. If $S$ is a sort, then ``$(\cdots, S^n, \cdots)$'' is shorthand for ``$(\cdots, S, \cdots, S, \cdots)$'' with $n$ copies of $S$.

If $f$ is a function, then $\dom(f)$ denotes the {\it  domain} $\{ x : \exists y \dt (\lrangle{x,y} \in f)\}$ of $f$, and $\mathrm{image}(f)$ denotes the {\it  image} $\{ y : \exists x \dt (\lrangle{x,y} \in f)\}$ of $f$. 

For the rest of the paper:
\begin{itemize}
\item Let $L$ be an arbitrary purely relational first-order language with only the one sort $\mathrm{Stat}$. Variables (e.g. $x, y, z, \vec{x}, \cdots$) of sort $\mathrm{Stat}$ are meant to refer to static individuals/tuples. This language is used throughout the paper to express properties of the static entities that do not change in the process being formalized. 
\item Let $L_\mathrm{Dyn}$ be a first-order language augmenting $L$ with 
\begin{itemize}
\item the sort $\mathrm{Dyn}_n$, for each $n < \omega$; 
\item the relation symbol $\leftleftarrows_n : (\mathrm{Stat}^n, \mathrm{Dyn}_n)$, for each $n < \omega$; 
\item possibly further relation symbols of any sort declaration.
\end{itemize}
Variables (e.g. $\xi, \upsilon, \zeta, \cdots$) of these sorts are meant to refer to dynamic tuples. $\vec{x} \leftleftarrows_n \xi$ is read ``$\xi$ {\it  manifests as} $\vec{x}$''.
\item Let $L_\mathrm{DynPred}$ be a first-order language augmenting $L_\mathrm{Dyn}$ with further relation symbols of any sort declaration.
\end{itemize}
The superscript on $\leftleftarrows$ is typically dropped when it is contextually clear or irrelevant.  The relation symbols $\leftleftarrows$ and the relation symbols of $L_\mathrm{DynPred} \setminus L_\mathrm{Dyn}$ are called {\it dynamic predicates}; the others are called {\it static predicates}. The dynamic predicates are used throughout the paper to express properties and relations that may change truth-value on an instance as the process unfolds, whereas the static predicates are used to express properties and relations that do not. Formally, in my framework of development structures (defined below in this section), the static predicates satisfy the stability axiom of convergent modal potentialism, as detailed in \S \ref{Sec:Background}.

The reason that $L$ is assumed to be purely relational is that function symbols carry an implicit commitment to the existence of a value, given any existing argument, which turns out to be cumbersome to manage in modal potentialist frameworks.\footnote{This is further explained in a footnote to the definition of the potentialist translation in \S \ref{Sec:Prel}.} This does not result in a significant loss of generality, since functions can be implemented as relations. In first-order logic, any function $f(\vec{x}, \vec{y})$ may be formalized as a relation $R(\vec{x}, \vec{y})$, along with the axiom $\forall \vec{x} \exists^! \vec{y} \dt R(\vec{x}, \vec{y})$, where the quantifier $\exists^! \vec{y}$ abbreviates a first-order formalization of ``there exists a unique tuple $\vec{y}$, such that''. If $f$ is a (partial) function, then $\dom(f)$ denotes the {\it  domain} $\{ x : \exists y \dt (\lrangle{x,y} \in f)\}$ of $f$, and $\mathrm{image}(f)$ denotes the {\it  image} $\{ y : \exists x \dt (\lrangle{x,y} \in f)\}$ of $f$. 

For the rest of the paper, let $\Lang$ be $L$, $L_\mathrm{Dyn}$, or $L_\mathrm{DynPred}$. $\Lang^\Diamond$ is the language used throughout the paper to express properties of dynamic entities. $\Lang_\textrm{StatPred}$ is the language of $\Lang$ restricted to static predicates. The choice of language $L$, $L_\mathrm{Dyn}$, or $L_\mathrm{DynPred}$ for $\Lang$ depends on which aspect of the process formalized are dynamic. This is made formally precise below by the definition of development model. But it may help the reader to explain their uses now: 
\begin{description}
	\item[$\Lang^\Diamond = L^\Diamond$] This cannot handle dynamic individuals, dynamic classes, or dynamic predicates. The only dynamics that this can handle is that the ambient structure increases through the process by passing to a superstructure. So once an element is introduced the truth value of an atomic formula does not change for that element as the process unfolds. This corresponds to Linnebo-style modal potentialism as developed in \cite{Lin13, LS19}.
	\item[$\Lang^\Diamond = L_\mathrm{Dyn}^\Diamond$] This can handle the above dynamics plus dynamic individuals or tuples, which can change their manifestation and thus the truth value of the atomic formulas they appear in as the process unfolds.
	\item[$\Lang^\Diamond = L_\mathrm{DynPred}^\Diamond$] This can handle the dynamics of $L_\mathrm{Dyn}$ plus dynamic predicates, whose truth value for a given static tuple can change as the process unfolds. 
\end{description}

Since function symbols carry an implicit commitment to the existence of a value, given any existing argument, they turn out to be cumbersome to manage in modal potentialist frameworks, which re-interpret existence as potential existence.\footnote{This is further explained in a footnote to the definition of the potentialist translation in \S \ref{Sec:Prel}.} I therefore follow the simple approach of side-stepping the issue by making the above assumption that $L$ is purely relational. This does not result in a significant loss of generality, since functions can be implemented as relations. In first-order logic, any function $f(\vec{x}, \vec{y})$ may be formalized as a relation $R(\vec{x}, \vec{y})$, along with the axiom $\forall \vec{x} \exists^! \vec{y} \dt R(\vec{x}, \vec{y})$, where the quantifier $\exists^! \vec{y}$ abbreviates a first-order formalization of ``there exists a unique tuple $\vec{y}$, such that''. Thanks to this fact, we will often use functional symbols in formulas, even though they are not strictly available in the language considered. In such cases, the formula is considered short-hand for the translation in terms of relation symbols, and the above axiom will be either an axiom of the theory we are working in, or satisfied by the model we are working in.

Let $u$ be a model of $\Lang$. $\mathrm{Th}(u)$ denotes the set of all $\Lang$-formulas $\phi$, such that $u \models \phi$. For any sort $S$ in its language, $S^u$ denotes the domain of $u$ of sort $S$. For a tuple $\vec{a}$, $\vec{a} \in u$ means that for each $1 \leq i \leq \lng(\vec{a})$, there is a sort $S$, such that $\pi_i(\vec{a}) \in \mathrm{S}^u$. For any $\Lang$-formula $\phi$, $\phi^u$ denotes its {\it  interpretation} $\{ \vec{a} \in u : u \models \phi(\vec{a})\}$ in $u$. Let $v$ be a model of $\Lang$. Let $K$ be a sublanguage of $\Lang$. 
\begin{itemize}
	\item $u$ is a $K$-{\it substructure} of $v$, denoted $u \subseteq_K v$ if the following hold:
	\begin{itemize}
		\item $S^u \subseteq S^v$, for all sorts $S$ in $\Lang$.
		\item For any sorts $S_1, \cdots, S_n$ and any relation symbol $R : (S_1, \cdots, S_n)$ in $K$, we have $R^u = R^v \cap (S_1^u \times \cdots \times S_n^u)$. 
	\end{itemize}
	\item $u$ is a {\it static substructure} of $v$ if $u \subseteq_{\Lang_\mathrm{StatPred}} v$.
	\item Suppose that $u, v$ are $L$-models. $u$ is a {\it $K$-elementary substructure of $v$}, denoted $u \preceq_K v$, if for all $\vec{a} \in u$ and all $\phi \in K$, $u \models \phi(\vec{a}) \Longleftrightarrow v \models \phi(\vec{a})$. 
\end{itemize}

In the definition of development model below it is required for a state $s'$ to be accessible from a state $s$ that the structure of $s$ be a static substructure of $s'$. This definition ensures that only the entities that are intended to be dynamic are so.

Let $W$ be a set of $\Lang$-models. Then $\bigcup W$ denotes the $\Lang$-model with $S^{\bigcup W} = \bigcup_{w \in W} S^w$ for each sort $S$ of $\Lang$, and with $R^{\bigcup W} = \bigcup_{w \in W} R^w$ for each relation symbol $R$.

Let $R \subseteq A \times A$ be a binary relation on a set $A$. 
\begin{itemize}
\item $R$ is {\it  reflexive} if $\forall x \in A \dt (x R x)$.
\item $R$ is {\it  antisymmetric} if $\forall x, y \in A \dt (x R y \wedge y R x \rightarrow x = y)$.
\item $R$ is {\it  transitive} if $\forall x, y, z \in A \dt \big( (x <^R y \wedge y <^R z) \rightarrow x <^R z \big)$.
\item $R$ is a {\it  partial order} if it is reflexive, antisymmetric and transitive.
\item $R$ is {\it  directed} if for all $a, a' \in A$, there is $b \in A$, such that $aRb$ and $a'Rb$.
\item $a \in A$ is a {\it  greatest element (with respect to $R$)} if $\forall x \in A \dt (x R a)$.
\item $R$ is {\it  unbounded} if $\forall x \in A \dt \exists y \in A \dt (x \neq y \wedge x R y)$.
\item $B \subseteq A$ is an {\it  up-set (with respect to $R$)} if $\forall b \in B \dt \forall a \in A \dt (b R a \rightarrow a \in B)$.
\item $B \subseteq A$ is a {\it  down-set (with respect to $R$)} if $\forall b \in B \dt \forall a \in A \dt (a R b \rightarrow a \in B)$.
\item $I \subseteq A$ is an {\it ideal (on $R$)} if $I$ is a down-set with respect to $R$ and $R \cap I \times I$ is directed. 
\item $D \subseteq A$ is {\it upward dense (with respect to $R$)} if $\forall a \in A \dt \exists d \in D \dt (a \leq d)$.
\item $I \subseteq A$ is a {\it generic ideal (on $R$)} if it is an ideal on $R$, such that for every upward dense $D$ with respect to $R$, we have $I \cap D \neq \varnothing$.
\end{itemize}

\begin{dfn}\label{Dfn:Dev_Mod}
A {\it  development structure/model} $\mathcal{M}$ ({\it  of} $\Lang^\Diamond$) consists of a non-empty set $\mathrm{states}^\mathcal{M}$ whose elements are called {\it  states} or {\it  worlds}, a relation $\leq^\mathcal{M}$ on $\mathrm{states}^\mathcal{M} \times \mathrm{states}^\mathcal{M}$ called the {\it  accessibility relation}, and a function $\mathrm{struct}^\mathcal{M}$ mapping each state $s$ to an $\Lang$-model $\mathrm{struct}^\mathcal{M}(s)$, such that:
\begin{enumerate}[{\rm (a)}]
\item $\leq^\mathcal{M}$ is reflexive, transitive, and directed;
\item $\forall n < \omega \dt \forall s \in \mathrm{states}^\mathcal{M} \dt \forall \delta \in (\mathrm{Dyn}_n)^{\mathrm{struct}^\mathcal{M}(s)} \dt \exists^! \vec{a} \in \mathrm{Stat}^{\mathrm{struct}^\mathcal{M}(s)} \dt \big( \langle \vec{a}, \delta \rangle \in (\leftleftarrows_n)^{\mathrm{struct}^\mathcal{M}(s)} \big)$;
\item For all $s, t \in \mathrm{states}^\mathcal{M}$, such that $s \leq^\mathcal{M} t$:
\begin{align*}
& \mathrm{struct}^\mathcal{M}(s) \statsub \mathrm{struct}^\mathcal{M}(t). \tag{$\statsub$-monotonicity}
\end{align*}
\end{enumerate}
\end{dfn}

Let $\mathcal{M}$ be a development model of $\Lang^\Diamond$. The structure $(\mathrm{states}^\mathcal{M}, \leq^\mathcal{M})$ is called a {\it  frame}. For any state $s$ of $\mathcal{M}$, $s$ also denotes $\mathrm{struct}^\mathcal{M}(s)$, as long as there is no risk for confusion.\footnote{The motivation for still defining development models with the $\mathrm{struct}^\mathcal{M}$ function is that there are cases where we need different states to map to the same structure.} Suppose that $S$ is an arbitrary sort of $\Lang$ other than $\mathrm{Stat}$. By $\subseteq_\mathrm{Dev}$-monotonicity, the map sending each state $s$ to the interpretation $S^{\mathrm{struct}^\mathcal{M}(s)}$ of $S$ at that state is constant. Therefore, we simply write $S^\mathcal{M}$ for that interpretation.

$\mathcal{M}$ is {\it finitary} if $\mathrm{Stat}^s$ is finite for all $s \in \mathrm{states}^\mathcal{M}$. Moreover, we define:
\begin{align*}
\mathrm{states}^\mathcal{M}_{\vec{a}} &=_\df \{s \in \mathrm{states}^\mathcal{M} : \vec{a} \in \mathrm{struct}^\mathcal{M}(s) \} \\
\bigcup \mathcal{M} &=_\df \bigcup \mathrm{image}(\mathrm{struct}^\mathcal{M})
\end{align*}
Let $u$ be an $L$-model. $\mathcal{M}$ is a {\it development model for $u$} if $u = \big(\bigcup \mathcal{M}) \restriction L$. Let $(A, \leq^A)$ be a reflexive, transitive, and directed relation on a set $A$. The {\em simple development for $u$ over $(A, \leq^A)$} is the development model $\mathcal{M}$ of $L^\Diamond$ defined as follows:
\begin{align*}
\mathrm{states}^\mathcal{M} =_\df A \\
\leq^\mathcal{M} =_\df {\leq}^A \\
\mathrm{struct}^\mathcal{M}(s) =_\df u
\end{align*}

Let $\kappa$ be a cardinal. $\mathcal{M}$ is {\it statewise} $\kappa$-ary if for each $s \in \mathrm{states}^\mathcal{M}$, we have $|\mathrm{Stat}^s| < \kappa$. $\mathcal{M}$ is {\it statewise finite} if it is statewise $\omega$-ary, and {\it statewise countable} if it is statewise $\omega_1$-ary. Note that for any first order model $u$ of $L$, there is a statewise finite development model for $u$: Indeed, the {\it full statewise finite development model for $u$}, denoted $\mathrm{FinDev}(u)$, is defined with the set of all finite subsets of the domain of $u$ as $\mathrm{states}^{\mathrm{FinDev}(u)}$, and with $\subseteq$ as $\leq^{\mathrm{FinDev}(u)}$, such that for any $s \in \mathrm{states}^{\mathrm{FinDev}(u)}$, $\mathrm{struct}^{\mathrm{FinDev}(u)}(s)$ is the $L$-substructure of $u$ on the domain $s$. It is immediate from Definition \ref{Dfn:Dev_Mod} that $\mathrm{FinDev}(u)$ is a well-defined development model.

Note that directedness and $\statsub$-monotonicity immediately yields the following proposition: 
\begin{prop}\label{Prop:Dev_Mod_Basics}
Let $\mathcal{M}$ be a development model of $\Lang$. If $\vec{a} \in \bigcup \mathcal{M}$, then $\mathrm{states}^\mathcal{M}_{\vec{a}}$ is a non-empty up-set with respect to $\leq^\mathcal{M}$.
\end{prop}

We define $S^\mathcal{M} =_\df \bigcup_{s \in \mathrm{states}^\mathcal{M}} S^{\mathrm{struct}(s)}$, for any state $s$, and any sort $S$.

A {\it static individual/tuple of $\mathcal{M}$} is a an individual/tuple of elements of $\mathrm{Stat}^{\bigcup \mathcal{M}}$. A {\it  dynamic ($n$-)tuple (of length $n$) of $\mathcal{M}$} is an element of $\mathrm{Dyn}_n^\mathcal{M}$. A {\it  dynamic individual} is a dynamic tuple of length $1$.

Let $s \in \mathrm{states}^\mathcal{M}$, let $n < \omega$, let $\delta \in \mathrm{Dyn}_n$, and let $\vec{a}$ be an $n$-tuple of elements of $\mathrm{Stat}^s$. $\delta$ {\it manifests as $\vec{a}$ in $s$} if $\langle \vec{a}, \delta \rangle \in (\leftleftarrows_n)^s$.

Let $s \in \mathrm{states}^\mathcal{M}$, let $\vec{a} \in \bigcup \mathcal{M}$ be an $n$-tuple, and let $\phi(\vec{x}) \in \Lang^\Diamond$. The satisfaction relation for development models is defined as follows, by structural induction: Let $\mathcal{M}$ be a development model of the language $\Lang^\Diamond$, let $s \in \mathrm{states}^\mathcal{M}$, let $n, m < \omega$, let $\vec{a}, \vec{b}$ be an $n$-tuple and $m$-tuple, respectively, of elements of arbitrary sorts of $\bigcup \mathcal{M}$, let $R$ be an $n$-ary relation symbol in $\Lang$, let $S$ be a sort of $\Lang$, and let $\phi(\vec{x}), \psi(\vec{y}) \in \Lang^\Diamond$, where $\vec{x}, \vec{y}$ are tuples of variables matching the length and sorts of $\vec{a}, \vec{b}$, respectively. 
\[
\begin{array}{lcl}
\mathcal{M}, s \models R(\vec{a}) &\Longleftrightarrow_\df& \vec{a} \in R^{\mathrm{struct}(s)}, \\
\mathcal{M}, s \models \neg \phi(\vec{a}) &\Longleftrightarrow_\df& s \not\models \phi(\vec{a}), \\
\mathcal{M}, s \models \phi(\vec{a}) \vee \psi(\vec{b}) &\Longleftrightarrow_\df& s \models \phi(\vec{a}) \vee s \models \psi(\vec{b}), \\
\mathcal{M}, s \models \exists x : S \dt \phi(x, \vec{a}) &\Longleftrightarrow_\df& \exists c \in S^{\mathrm{struct}(s)} \dt s \models \phi(c, \vec{a}), \\
\mathcal{M}, s \models \Diamond \phi(\vec{a}) &\Longleftrightarrow_\df& \exists t \dt \big( s \leq^\mathcal{M} t \wedge t \models \phi(\vec{a}) \big).
\end{array}
\]
When the context is clear we write ``$s \models \cdots$'' for ``$\mathcal{M}, s \models \cdots$''. $\mathcal{M} \models \phi(\vec{a})$ is defined as $\forall u \in \mathrm{states}^\mathcal{M} \dt \big( u \models \phi(\vec{a}) \big)$. Note that if $\vec{a}$ is a tuple of elements of appropriate sorts of $\mathrm{struct}^\mathcal{M}(s)$ and $\phi(\vec{x}) \in \Lang$, then $s \models \phi(\vec{a}) \Longleftrightarrow \mathrm{struct}(s) \models \phi(\vec{a})$.

By definition of development models, for each state $s$ of $\mathcal{M}$, each dynamic tuple $\delta$ of $\mathcal{M}$ manifests as a unique tuple $\vec{a}$ in $s$. This is meant to formalize the philosophical nature of dynamic entities explained in \S \ref{Sec:Char_Dev_Proc}. In particular, the type coherence of dynamic entities leads to a design requirement that we should be able to regard $\delta$ as having a type which is the dynamic version of the type of the tuples it manifests as in the states. Therefore, any formula $\phi(\vec{x})$, where $\vec{x} : (\mathrm{Stat}^n)$, should not only be applicable to static $n$-tuples of $\mathcal{M}$, but also to dynamic $n$-tuples of $\mathcal{M}$. On the other hand, from a strictly formal perspective, we want to keep the integrity of the sort declarations intact, not allowing variables to range over multiple sorts. So to achieve the above purpose, we introduce a translation: Let $\phi(\vec{x})$ be an $\Lang$-formula, where $\vec{x} : (\mathrm{Stat}^n)$. Without loss of generality, we may assume that no occurrence in $\phi$ of a variable from $\vec{x}$ is bound. The {\it dynamic substitution of $\xi : \mathrm{Dyn}_n$ for $\vec{x}$ in $\phi(\vec{x})$}, denoted $\phi(\vec{x} \tilde{/} \xi)$ (or $\phi(\xi)$ when the particular variables dynamically substituted are clear), is the formula obtained from $\phi(\vec{x})$ by replacing each atomic subformula $\psi$ of $\phi$ by the formula $\exists \vec{x} \dt \big( \vec{x} \leftleftarrows_n \xi \wedge \psi \big)$.\footnote{Formally this is defined in the usual way by structural induction, where the base case is as stated, and the cases for the connectives, quantifiers, and operators are trivial.} 

\begin{prop}
Let $\mathcal{M}$ be a development model of $L_{\mathrm{Dyn}}$, let $s \in \mathrm{states}^\mathcal{M}$, and let $\phi(\vec{x}) \in \Lang^\Diamond$, where $\vec{x} : (\mathrm{Stat}^n)$ do not occur bound in $\phi$. Moreover, let $\vec{a} \in \mathrm{Stat}^s$ be an $n$-tuple, and let $\delta \in (\mathrm{Dyn}_n)^s$ that manifests as $\vec{a}$ in $s$. Then:
\begin{align*}
\mathcal{M}, s \models \phi(\vec{a}) \Longleftrightarrow \mathcal{M}, s \models \phi(\vec{a} \tilde{/} \delta)
\end{align*}
\end{prop}
\begin{proof}
The only non-trivial case is the base case, and this case follows from that $\delta$ manifests as $\vec{a}$ in $s$.
\end{proof}

Any development model of $\Lang^\Diamond$ converges to its union, with respect to $L$, in the following sense:

\begin{prop}\label{Prop:Converge}
Suppose that $\mathcal{M}$ is a development model of $\Lang^\Diamond$. Let $\phi$ be a literal in $\Lang_\mathrm{StatPred}$, and let $\vec{a} \in \bigcup \mathcal{M}$ matching the free variables of $\phi$. The following are equivalent:
\begin{enumerate}[{\rm (a)}]
\item\label{Prop:Converge_union} $\bigcup \mathcal{M} \models \phi(\vec{a})$ 
\item\label{Prop:Converge_exists} $\exists s \in \mathrm{states}^\mathcal{M} \dt \big( \mathrm{struct}^\mathcal{M}(s) \models \phi(\vec{a}) \big)$
\item\label{Prop:Converge_forall} $\mathcal{M} \models \Diamond \Box \phi(\vec{a})$
\end{enumerate}
\end{prop}
\begin{proof}
By directedness and $\statsub$-monotonicity, we have (\ref{Prop:Converge_exists}) $\Longleftrightarrow$ (\ref{Prop:Converge_forall}).

Suppose that $\phi$ is an atomic $\Lang_\mathrm{StatPred}$-formula $R(x_1, \cdots, x_n)$. By definition of $\bigcup \mathcal{M}$:
\begin{align*}
\bigcup \mathcal{M} \models R(\vec{a}) &\Longleftrightarrow \exists s \in \mathrm{states}^\mathcal{M} \dt \big( \mathrm{struct}^\mathcal{M}(s) \models R(\vec{a}) \big) 
\end{align*}

Suppose, on the other hand, that $\phi$ is a negated atomic $\Lang_\mathrm{StatPred}$-formula $\neg R(x_1, \cdots, x_n)$. By definition of $\bigcup \mathcal{M}$:
\begin{align*}
\bigcup \mathcal{M} \models \neg R(\vec{a}) &\Longleftrightarrow \forall s \in \mathrm{states}^\mathcal{M}_{\vec{a}} \dt \big( \mathrm{struct}^\mathcal{M}(s) \models \neg R(\vec{a}) \big)
\end{align*}
By directedness and $\statsub$-monotonicity there is a state $t$ containing $\vec{a}$. So 
\begin{align*} 
\forall s \in \mathrm{states}^\mathcal{M}_{\vec{a}} \dt \big( \mathrm{struct}^\mathcal{M}(s) \models \neg R(\vec{a}) \big) \Longrightarrow \mathrm{struct}^\mathcal{M}(t) \models \neg R(\vec{a}).
\end{align*}
Conversely, by directedness and $\statsub$-monotonicity: 
\begin{align*}
\exists s \in \mathrm{states}^\mathcal{M}_{\vec{a}} \dt \big( \mathrm{struct}^\mathcal{M}(s) \models \neg R(\vec{a}) \big) \Longrightarrow \bigcup \mathcal{M} \models \neg R(\vec{a})
\end{align*}
Thus, we have (\ref{Prop:Converge_union}) $\Longleftrightarrow$ (\ref{Prop:Converge_exists}).
\end{proof}

Although this paper proposes an alternative to the usual functional representation of entities developing through a process, the functional representation is useful in the meta-theory: Let $n < \omega$. A {\it  functional dynamic ($n$-)tuple (of length $n$) of $\mathcal{M}$} is a function $D$, such that 
\begin{itemize}
\item $\dom(D) = \mathrm{states}^\mathcal{M}$,
\item $\forall s \in \dom(D) \dt \big( D(s) \in (\mathrm{Stat}^s)^n \big)$.
\end{itemize}
A {\it  functional dynamic individual} is a dynamic tuple of length $1$.

Note that every dynamic tuple $\delta$ of $\mathcal{M}$ induces a functional dynamic tuple $D_\delta$ mapping each state $s$ to the unique tuple $\vec{a}$, such that $\delta$ manifests as $\vec{a}$ in $s$.

Let $\mathcal{D}$ be a set of functional dynamic tuples (possibly of various lengths) of $\mathcal{M}$. Define {\it  the least development model extending $\mathcal{M}$ by $\mathcal{D}$}, denoted $\tilde{\mathcal{M}}(\mathcal{D})$, as the extension of $\mathcal{M}$ obtained by adding an element $\delta_{D}$ of sort $\mathrm{Dyn}_{\lng(D)}$, injectively for each $D \in \mathcal{D}$, along with the following extension of the interpretation of $\leftleftarrows_n$, for each $n < \omega$, for each $D \in \mathcal{D}$ of length $n$, and for each $s \in \mathrm{states}^\mathcal{M}$:
\begin{align*}
\langle \vec{a}, \delta_D \rangle \in (\leftleftarrows_n)^{\mathrm{struct}^{\tilde{\mathcal{M}}(\mathcal{D})}(s)} \Longleftrightarrow_\df D(s) = \vec{a}
\end{align*}
$\tilde{\mathcal{M}}$ denotes the least development model extending $\mathcal{M}$ by the set of all functional dynamic individuals, and $\tilde{\mathcal{M}}_{< \omega}$ denotes the least development model extending $\mathcal{M}$ by the set of all functional dynamic tuples.

Although not formally necessary, it is convenient to define for all states $s$, functional dynamic tuples $D \in \mathcal{D}$, and formulas $\phi(\vec{x}) \in \Lang^\Diamond$:
\begin{align*}
\mathcal{M}, s \models \phi(D) &\Longleftrightarrow_\df \tilde{\mathcal{M}}(\mathcal{D}), s \models \phi(\delta_D) 
\end{align*}

\section{Background}\label{Sec:Background}

This section addresses and leverages previous research that is related to the work of this paper. First of all, the formalization of dynamic objects is analogous to Carnap's formalization of individual concepts in \cite{Car56}. Carnap's semantics satisfies $\mathsf{S5}$, and so the accessibility relation is an equivalence relation, in contrast to the directed partial order employed for development models. We can thus think of the notion of dynamic objects as a procedural cousin of the notion of individual concepts: intensional in the modal sense, but extensional in their directed development (this extensionality is typically captured by an equivalence relation in the applications throughout this paper). This section is primarily concerned with the highly relevant potentialist view of infinity and the foundations of mathematics. In particular, the mirroring theorems of \cite{LS19, HL19} are generalized to the framework of development models. Moreover, there are some affinities between this framework, and the frameworks developed in \cite{HS19, Kri19, BLS22, Bra23} which are explained below.

The potentialist view of infinity originates with Aristotle, who formulated it in the {\it  Physics} ca. 350 B.C. \cite[Book 3, Part 6]{A}: ``It is not what has nothing outside it that is infinite, but what always has something outside it. [\ldots] This condition alone is not sufficient: it is necessary also that the next part which is taken should never be the same. [\ldots] Our definition then is as follows: A quantity is infinite if it is such that we can always take a part outside what has been already taken.''

A modal framework formalizing a potentialist position for arithmetic is developed in \cite{LS19}, along the lines of an analogous development for Zermelo-Fraenkel ($\mathsf{ZF}$) style set theory introduced in \cite{Lin13}. The general modal theory of this framework is $\mathsf{S4.2}$ over first-order logic plus the axiom schema of $\textsf{Stability}$ asserting $\phi \rightarrow \Box \phi$, for each literal $\phi$. I call this theory $\textsf{Convergent Modal Potentialism}$ ($\mathsf{CMP}$). Recall that $\mathsf{S4.2}$ is $\mathsf{S4}$ plus the axiom schema $\mathsf{G}$ asserting $\Diamond\Box \phi \rightarrow \Box\Diamond \phi$, for all $\phi$. $\mathsf{CMP}$ proves the schema $\textsf{Converse Barcan Formula}$ ($\mathsf{CBF}$) $\exists x \Diamond \phi \rightarrow \Diamond \exists x \phi$, for all $\phi$. There is a direct bridge between first-order logic and $\mathsf{CMP}$: The {\it  potentialist translation} is a mapping that sends each first-order formula $\phi$ to its {\it  full modalization}, denoted $\phi^\Diamond$, where $\phi^\Diamond$ is the formula obtained by replacing each occurrence of the first-order quantifiers $\exists$ and $\forall$ in $\phi$ by $\Diamond\exists$ and $\Box\forall$, respectively (note that the second-order quantifiers remain unchanged).\footnote{If first-order function symbols were admitted in the languages, then this translation would need to be modified to a more complicated one. For example, consider the axiom $\forall x \mathrm{S}(x) \neq x$ of $\PA$, with `$\mathrm{S}$' as function symbol. When $\mathsf{CMP}$ is extended to a modal theory of arithmetic, the translation $\Box\forall x \mathrm{S}(x) \neq x$ should not be provable, because it is naturally interpreted as committed to that every number has a successor, contradicting the potentialist view to be formalized. However, the translation of its relational counterpart is $\Box\forall x \forall y \big( \mathrm{S}(x, y) \rightarrow x \neq y \big)$, which should be provable in the potentialist framework.} This notation is extended to any first-order $L$-theory $T$, so that $T^\Diamond$ denotes the $L^\Diamond$-theory extending $\mathsf{CMP}$ by $\phi^\Diamond$, for each $\phi$ in $T$. The fundamental theorem of the $\mathsf{CMP}$ framework is the Mirroring Theorem \cite[Theorem 1]{LS19}, here expressed as follows:

\begin{thm}[Syntactic Mirroring Theorem]\label{Thm:Synt_Mirror}
For any $L$-theory $T$, and any $\phi \in L$:
\[
T \vdash \phi \iff T^\Diamond \vdash \phi^\Diamond
\]
\end{thm}

The framework is endowed with the usual possible worlds semantics, to which the static part of any development model belongs (see Corollary \ref{Cor:Dev_Mod_CMP}) below.

Moreover, we immediately obtain the following semantic version of the Mirroring Theorem from \cite[Theorem 1]{HL19}:

\begin{thm}[Semantic Mirroring Theorem]\label{Thm:Sem_Mirror}
For any development model $\mathcal{M}$ of $\Lang^\Diamond$, any $\vec{a} \in \bigcup \mathcal{M}$, and any $\phi \in \Lang_\mathrm{StatPred}$: \[
\bigcup \mathcal{M} \models \phi(\vec{a}) \iff \mathcal{M} \models \phi^\Diamond(\vec{a})
\]
\end{thm}
\begin{proof}
	This follows from \cite[Theorem 1]{HL19}. We just need to show that $\mathcal{M} \restriction \Lang_\mathrm{StatPred}^\Diamond$ is a {\it potentialist account} of $\bigcup \big( \mathcal{M} \restriction \Lang_\mathrm{StatPred}^\Diamond \big)$ in the sense of the definition in \cite[\S 1]{HL19}. That is, we need to show that for any state $s$ in $\mathrm{states}(\mathcal{M})$: 
	\begin{enumerate}
		\item $\mathrm{Struct}^\mathcal{M}(s) \restriction \Lang_\mathrm{StatPred}$ is a substructure of $\bigcup \big( \mathcal{M} \restriction \Lang_\mathrm{StatPred}^\Diamond \big)$.
		\item For every $\vec{b} \in \bigcup \big( \mathcal{M} \restriction \Lang_\mathrm{StatPred}^\Diamond \big)$, there is a state $s' \geq^\mathcal{M} s$, such that $\vec{b} \in \mathrm{Struct}^\mathcal{M}(s')$.
	\end{enumerate}
	The first follows from the definition of $\bigcup \big( \mathcal{M} \restriction \Lang_\mathrm{StatPred}^\Diamond \big)$. The second follows from the definition of $\bigcup \big( \mathcal{M} \restriction \Lang_\mathrm{StatPred}^\Diamond \big)$, and from directedness and item (c) of Definition \ref{Dfn:Dev_Mod}.
\end{proof}

\begin{cor}\label{Cor:Dev_Mod_CMP}
For each development model $\mathcal{M}$ of $\Lang^\Diamond$:
\[
\mathcal{M} \models \mathsf{CMP}\restriction \Lang_\textrm{StatPred}^\Diamond,
\]
where $\mathsf{CMP}\restriction \Lang_\textrm{StatPred}^\Diamond$ is the theory $\mathsf{CMP}$ expressed in the language $\Lang_\textrm{StatPred}^\Diamond$ (i.e. with its schemata ranging over formulas of that language).	
\end{cor}
\begin{proof}
	Immediate from Theorems \ref{Thm:Synt_Mirror} and \ref{Thm:Sem_Mirror}.´
\end{proof}

Let $u$ be any $L$-model. By Theorem \ref{Thm:Sem_Mirror}, $\mathrm{FinDev}(u) \models \mathrm{Th}(u)^\Diamond$. 

Now take $L$ to be the language of Peano Arithmetic ($\mathsf{PA}$) formalized in the language with a unary predicate $\mathrm{Z}(x)$, expressing that $x = 0$, a binary predicate $\mathrm{S}(x, y)$, expressing that $y$ is the successor of $x$, and ternary predicates $+$, $\times$, formalizing the operations of addition, and multiplication, respectively. The {\it  basic development model of arithmetic}, is the $L^\Diamond$-structure $\mathcal{M}$ on the frame $(\mathbb{N}, \leq)$, where $\mathrm{struct}(n)$ is the substructure of $\mathbb{N}$ with domain $\{0, \cdots, n\}$. It follows from Proposition \ref{Thm:Sem_Mirror} that $\mathcal{M}$ models $\PA^\Diamond$. Moreover, $\mathcal{M} \models \Box \forall x \dt \Diamond \exists y \dt S(x, y)$, and $\mathcal{M} \models \Box \exists y \dt \forall x \dt x \leq y$. These may be viewed as the twin Aristotelian claims that infinity is potential, and that the actual is necessarily finite, respectively. The modal potentialist account is developed in \cite{LS19}.

Let us consider how the $\mathsf{CMP}$ framework balances between finitism and infinitism, and combines some of their respective advantages. The property of the above development model that every world is finite, but that there are infinitely many worlds, reflects the ``classical finitist'' position, but conflicts with the ``strict finitist'' position, in the terminology of \cite[pp. 6--8]{Til89}, because the classical finitist accepts the existence of all finite structures, but no infinite structures. There are many philosophical options for how to construe the potentialist modality of the framework. In order to remain neutral on metaphysical debates, it is convenient to construe it as an intepretational modality. The position is then that only finite interpretations of arithmetic are admitted. One view motivating this is that only finite interpretations of arithmetic are meaningful; another is that only finite structures are metaphysically possible. These positions are finitist in character. On the other hand, the Mirroring Theorems accommodate elements of infinitism: They reveal direct syntactic and semantic correspondances between the $\textsf{CMP}$ framework and the usual first-order framework of mathematics. In particular, they show that $\PA^\Diamond$ has all the ``reasoning power'' of $\PA$, and that the basic model of potentialist arithmetic has all the ``satisfaction power'' of the standard model of arithmetic, respectively.

In contrast to the above convergent modal potentialism, \cite{BLS22, Bra23} develop a formal framework for divergent modal potentialism, applied to formalize free choice sequences (a concept of Brouwerian intuitionism). The idea of such choice sequences is that they are necessarily of finite length, but can grow by freely choosing a next value to append. Accordingly, the formalization of \cite{BLS22, Bra23} allows the value of such a choice sequence to diverge, when passing from a possible world to different accessible worlds. As the identity of the choice sequence itself remains the same, this behavior requires that the underlying frame be a non-directed partial order. Moreover, it involves a modal operator of inevitability, deemed to be required by divergent potentialism. Their formalization is a natural adaptation of the $\mathsf{CMP}$-formalization to divergent potentialism. \cite{Kri19} outlines a related approach to free choice sequences, based on Beth models. Such models are again not assumed to be directed.
The above models are similar to development models in that an entity is formalized which can take on different values in different worlds, but they differ fundamentally in that the former are non-directed and the latter are directed. Directedness is a key property of development processes that enables the above mirroring theorems for classical logic,\footnote{\cite{BLS22} proves mirroring theorems for intuitionistic logic by means of a certain translation to their framework.} and many of the results of this paper. In particular, directedness is essential to formalizing the teleological convergence that many development processes have.

Another notion in the philosophy of mathematics which concerns entities that can take on a variety of values is that of arbitrary objects. Such a notion is developed with philosophical rigor in \cite{Fin85a, Fin85b}. This has inspired formal work, including \cite{HS19, VY24}, which employ modal logic to formalize arbitrary objects. In particular, Horsten and Speranski \cite{HS19} treat arbitrary numbers as Carnapian individual concepts over the space of all permutations of $\mathbb{N}$, yielding an $\mathsf{S5}$-modal framework. An arbitrary object is thus quite different in character from a dynamic object, in that arbitrary objects have no directed development connecting the various values it can take, while dynamic objects are procedural entities, formalized as unfolding along a directed frame.

\section{Teleological development in general}\label{Sec:Tel_Dev}

Typically, one is interested in teleological conditions that a process develops in order to meet. In line with the blueprint given in \S \ref{Sec:Char_Dev_Proc}, we formalize this notion as follows: Let $\mathcal{M}$ be a development model of $\Lang^\Diamond$, let $\vec{a} \in \mathcal{M}$, and let $\phi \in \Lang^\Diamond$. {\it $\mathcal{M}, \vec{a}$ satisfy $\phi$ teleologically}, denoted $\mathcal{M} \models^\mathrm{tele} \phi(\vec{a})$, if $\mathcal{M} \models \Diamond \Box \phi(\vec{a})$. Let $K$ be a sublanguage of $\Lang^\Diamond$. The $K$-{\it teleology of $\vec{a}$ over $\mathcal{M}$} is the set of formulas $\phi$ in $K$, such that $\mathcal{M} \models^\mathrm{tele} \phi(\vec{a})$. (If ``$K$-'' is omitted it means that $K = \Lang^\Diamond$. If ``of $\vec{a}$'' is omitted it means that $\vec{a} = \langle\rangle$.)

As an example, let $\mathrm{Lit}$ be the set of literals of $L$ and let $u$ be an $L$-model. By Proposition \ref{Prop:Converge}, there is a statewise finite development model for $u$, namely $\mathrm{FinDev}(u)$, such that for any $\vec{a} \in u$, the $\mathrm{Lit}$-teleology of $\vec{a}$ over $\mathrm{FinDev}(u)$ is the set of $\phi \in \mathrm{Lit}$, such that $u \models \phi(\vec{a})$. This raises:
\begin{que*}
Let $u$ be an $L$-model. Is there a statewise ``reasonably small'' development model $\mathcal{M}$ for $u$, such that for any $\vec{a} \in u$, the $L$-teleology of $\vec{a}$ over $\mathcal{M}$ is the set of $\phi \in L$, such that $u \models \phi(\vec{a})$?
\end{que*}

The following refinement of the Downward L\"owenheim--Skolem theorem gives a strongly coherent affirmative answer to this question:

\begin{thm}[L\"owenheim--Skolem-style]
Let $u$ be a model of $L$. There is a statewise countable development model $\mathcal{M}$ for $u$, such that $\mathcal{M} \models^\mathrm{tele} \mathrm{Th}(u)$. Moreover:
\begin{align*}
& \mathrm{states}^\mathcal{M} \textrm{ is the set of all countable sublanguages of $L$.} \\
& \forall s, t \in \mathrm{states}^\mathcal{M} \dt \big( s \leq^\mathcal{M} t \leftrightarrow s \subseteq t \big) \\
& \forall s \leq^\mathcal{M} t \in \mathrm{states}^\mathcal{M} \dt \big( \mathrm{struct}^\mathcal{M}(s) \preceq_s \mathrm{struct}^\mathcal{M}(t) \preceq_t u \big) \tag{$\preceq$-monotonicity}
\end{align*}
\end{thm}
\begin{proof}
This is proved by the well-known technique of Skolem hulls, introduced e.g. in \cite[Ch. 3]{Hod93}. We fix a skolemization $L^*, T^*$ of $L, \mathrm{Th}(u)$, and let $v$ be an expansion of $u$ modeling $T^*$ (using Theorem 3.1.2 and Corollary 3.1.4 in \cite{Hod93}). $\mathcal{M}$ is constructed as follows: Let $\mathrm{states}^\mathcal{M}$ be the set of all countable sublanguages of $L$, and let $\leq^\mathcal{M}$ be the subset-relation. For each $s \in \mathrm{states}^\mathcal{M}$, we define $s^* \subseteq L$ as the skolemization of $s$. Moreover, for each $s \in \mathrm{states}^\mathcal{M}$, we define $\mathrm{struct}^\mathcal{M}(s)$ as the following $L$-structure: it interprets $s$ as the $s$-reduct of the Skolem hull of $\varnothing$ in the $s^*$-reduct of $v$, and it interprets each relation symbol in $L \setminus s$ by $\varnothing$. For each state $s$, since $s$ is a countable language, $s^*$ is a countable language and $\mathrm{struct}^\mathcal{M}(s)$ is a countable structure. Let $s \leq^\mathcal{M} t \in \mathrm{states}^\mathcal{M}$. By definition, we have that $\mathrm{struct}^\mathcal{M}(s)$ is the $s$-reduct of the Skolem hull of $\varnothing$ in the $s^*$-reduct of $t$. So by Theorem 3.1.1 in \cite{Hod93}, we obtain $\mathrm{struct}^\mathcal{M}(s) \preceq_s \mathrm{struct}^\mathcal{M}(t) \preceq_t u$, establishing $\prec$-monotonicity. 

To see that $\mathcal{M} \models^\mathrm{tele} \mathrm{Th}(u)$, let $\phi \in L$ and let $s \in \mathrm{states}^\mathcal{M}$. Then there is $s' \geq^\mathcal{M} s$, such that $\phi \in s'$. By $\preceq$-monotonicity, we have for all $t  \geq^\mathcal{M} s'$ that $s' \models \phi \Longleftrightarrow t \models \phi \Longleftrightarrow u \models \phi$. Thus, $\mathcal{M} \models \Diamond \Box \phi \Longleftrightarrow u \models \phi$.
\end{proof}

The above theorem is optimal in the following senses: If $u$ is finite, then the question is trivially answered by the development model with a single state mapped to the model $u$. If $u$ is infinite, then there is no finite model $v$, such that $v \preceq_\mathrm{FOL} u$, where $\mathrm{FOL}$ is the basic language of first order logic. So there is no chance to improve the theorem to ``statewise finite'' while maintaining $\preceq$-monotonicity. Even if $\preceq$-monotonicity is dropped, the formula $\exists x \forall y (y \leq x)$ shows that there is no statewise finite development model $\mathcal{M}$ for $\mathbb{N}$, such that $\mathcal{M} \models^\mathrm{tele} \mathrm{Th}(\mathbb{N})$. So for arbitrary $u$, there is no chance to improve the theorem to ``statewise finite''. However, if $u$ is just an infinite model of $\mathrm{FOL}$ (i.e. with no structure), then it is clear from basic first order model theory that there is a statewise finite development model for $u$ that teleologically satifies the theory of $u$: indeed, $\mathrm{FinDev}(u) \models^\mathrm{tele} \mathrm{Th}(u)$.  This raises:
\begin{que*}
Is there a natural and general class of $L$-theories $T$, and infinite models $u \models T$, such that there is a statewise finite development model for $u$ that $L$-teleologically satifies $T$?
\end{que*}

Theorem \ref{Thm:Statewise_finite_dev_mod} below gives an affirmative answer in terms of logical complexity. A theory $T$ is a {\it $\Sigma_2$-theory} if each $\phi \in T$ is of the form $\exists x_1 \cdots \exists x_m \forall y_1 \cdots \forall y_n \dt B_\phi$, where $B_\phi$ is a quantifier free formula.

\begin{thm}\label{Thm:Statewise_finite_dev_mod}
Let $u$ be an infinite model of a $\Sigma_2$-theory $T$. Then $\mathrm{FinDev}(u)$ is a statewise finite development model for $u$ that $L$-teleologically satifies $T$.
\end{thm}
\begin{proof}
Let $\phi \in T$, and let $s$ be a state of $\mathrm{FinDev}(u)$ (recall that $s$ is a finite subset of the domain of $u$). Then $\phi$ is of the form $\exists x_1 \cdots \exists x_m \forall y_1 \cdots \forall y_n \dt B_\phi(x_1, \cdots, x_m, y_1, \cdots, y_n)$, where $B_\phi$ is a quantifier free formula. Since $u \models T$, there are $a_1, \cdots, a_m \in u$, such that $u \models \forall y_1 \cdots \forall y_n \dt B_\phi(a_1, \cdots, a_m, y_1, \cdots, y_n)$. Let $t$ be the finite subset $s \cup \{ a_1, \cdots, a_m \}$ of the domain of $u$. By the form of $\phi$, we have for any state $t' \geq t$ of $\mathrm{FinDev(u)}$ that $\mathrm{struct}^{\mathrm{FinDev(u)}}(t') \models \phi$. Thus, $\mathrm{FinDev}(u) \models^{\mathrm{tele}} \phi$, as desired. 
\end{proof}

\section{Reals as dynamic rationals}\label{Sec:Reals}

Assume that $L$ is the language of a relational axiomatization $\mathsf{DLOF}$ of the first-order theory of dense linearly ordered fields (i.e. with addition and multiplication formalized as ternary relations). Let $\mathbb{Q}$ be the field of the rational numbers. $\mathbb{Q}$ is well-known to be the prime model of $\mathsf{DLOF}$.  Let $\mathbb{Q}_+$ be the set of positive rationals. Let $\mathbb{R}$ be the relational structure of the field of real numbers. The next theorem shows that for any finitary development model $\mathcal{Q}$ of $L^\Diamond$ for $\mathbb{Q}$, we have that $\tilde{\mathcal{Q}}$ defines $\mathbb{R}$ in a certain precise sense. The proof is analogous to the well-known construction of $\mathbb{R}$ as the set of all Cauchy convergent sequences of rationals modulo Cauchy equivalence. However, in the present setting it is more natural to use the generalization of sequences to nets. Nets are discussed in introductory topology texts, e.g. \cite[p. 187]{Mun00}. A {\it net} in a topological space is a function from a directed partial order into the space. Thus, any functional dynamic individual of a development model for $\mathbb{Q}$ is a net. 

Let $\mathcal{T}$ be the standard topological space of $\mathbb{Q}$ or $\mathbb{R}$. Let $\langle P, \leq^P \rangle$ be a directed partial order, and let $N, N', N'' : P \rightarrow \mathcal{T}$ be nets. In this context: $N$ is {\it Cauchy convergent} if $\forall z \in \mathbb{Q}_+ \dt \exists p \in P \dt \forall p', p'' \geq^P p \dt \big( |N(p') - N(p'')| < z \big)$. $N$ and $N'$ are {\it Cauchy equivalent}, denoted $N \sim N'$, if $\forall z \in \mathbb{Q}_+ \dt \exists p \in P \dt \forall p' \geq^P p \dt \big( |N(p') - N'(p')| < z \big)$. $N$ {\it converges to} $a \in \mathcal{T}$ if $\forall z \in \mathbb{Q}_+ \dt \exists p \in P \dt \forall p' \geq^P p \dt \big( |a - N(p')| < z \big)$. Addition and multiplication of nets are defined component-wise by $+(N, N', N'') \leftrightarrow_\df \forall p \in P \dt +(N(p), N'(p), N''(p))$ and $\times(N, N', N'') \leftrightarrow_\df \forall p \in P \dt \times(N(p), N'(p), N''(p))$. The following facts are easy to see (and well-known): Cauchy equivalence is an equivalence relation. $N$ is Cauchy convergent iff it converges to a number in $\mathbb{R}$. If $N$ converges to a number in $\mathcal{T}$, then it converges to exactly one number in $\mathcal{T}$, denoted $r_N$. If $N' \sim N$, then $N$ converges to $a$ in $\mathcal{T}$ iff $N'$ converges to $a$ in $\mathcal{T}$. If $M, M', M''$ are nets that are Cauchy equivalent to $N, N', N''$, respectively, then $+(M, M', M'')$ iff $+(N, N', N'')$, and $\times(M, M', M'')$ iff $\times(N, N', N'')$.

\begin{lemma}\label{Lem:Nets_R}
Let $\langle P, \leq^P \rangle$ be a directed partial order. Let $\mathcal{N}$ be a set of Cauchy convergent nets $N : P \rightarrow \mathbb{Q}$, such that for each $r \in \mathbb{R}$, there is $N \in \mathcal{N}$, such that $r_N = r$. Then:
\begin{align*}
\langle \mathcal{N}, +, \times \rangle / {\sim} \cong \langle \mathbb{R}, +, \times \rangle
\end{align*}
\end{lemma}
\begin{proof}
By the above facts, $\langle \mathcal{N}, +, \times \rangle / {\sim}$ is well-defined. The proposed isomorphism is induced by $N \mapsto r_N$, for each $N \in \mathcal{N}$. By the facts cited above this lemma, this induces a well-defined map from $\mathcal{N} / \sim$ to $\mathbb{R}$, which is injective. By the assumption of the theorem, this map is also surjective. Finally, by the facts cited above this lemma, it preserves all relevant structure. Thus, it is an isomorphism.
\end{proof}

\begin{thm}
Let $\mathcal{Q}$ be any finitary development model of $L^\Diamond$ for $\mathbb{Q}$. Then $\tilde{\mathcal{Q}}$ defines $\mathbb{R}$, in the following sense: There are formulas $\rho(\xi)$, $\xi \sim \upsilon$, $\tilde{+}(\xi, \upsilon, \zeta)$, $\tilde{\times}(\xi, \upsilon, \zeta)$ in $L_\mathrm{Dyn}^\Diamond$, such that 
\begin{align*}
\langle \rho^{\tilde{\mathcal{Q}}}, \tilde{+}^{\tilde{\mathcal{Q}}}, \tilde{\times}^{\tilde{\mathcal{Q}}} \rangle  / {\sim^{\tilde{\mathcal{Q}}}} \cong \mathbb{R}.
\end{align*}
\end{thm}
\begin{proof}
Let $\mathcal{D}$ be the set of all functional dynamic individuals over $\mathcal{Q}$. Recall that $\tilde{\mathcal{Q}} = \tilde{\mathcal{Q}}(\mathcal{D})$, and that every dynamic individual $D \in \mathcal{D}$ is a net in $\mathbb{Q}$ and in $\mathbb{R}$.

The key insight is that natural notions of Cauchy convergence and Cauchy equivalence are definable in $\tilde{\mathcal{Q}}$, by the following formulas:
\begin{align*}
\rho(\xi) &\leftrightarrow_\df \Box \forall z > 0 \dt \Diamond \Box \exists x \dt \Big( x \leftleftarrows \xi \wedge \Box \big( (|x - \xi| < z)^\Diamond \big) \Big) \\
\xi \sim \upsilon &\leftrightarrow_\df \Box \forall z > 0 \dt \Diamond \Box \big( (|\xi - \upsilon| < z)^\Diamond \big)
\end{align*}
The formula $(|x - y| < z)^\Diamond$ is obtained by the potentialist translation (see \S \ref{Sec:Background}) from the formula $|x - y| < z$, which in turn is short-hand for an $L$-formula with quantifiers. Thus, it is worth noting that $(|x - y| < z)^\Diamond$ involves modal operators. The formulas $(|x - \xi| < z)^\Diamond$ and $(|\xi - \upsilon| < z)^\Diamond$ are obtained from $(|x - y| < z)^\Diamond$ by dynamic substitution (see \S \ref{Sec:Prel}). 

Addition and multiplication of dynamic individuals are defined as follows:
\begin{align*}
\tilde{+}(\xi, \upsilon, \zeta) &\leftrightarrow_\df \Box \mathop{+}(\xi, \upsilon, \zeta) \\
\tilde{\times}(\xi, \upsilon, \zeta) &\leftrightarrow_\df \Box \mathop{\times}(\xi, \upsilon, \zeta)
\end{align*}

Let $\mathcal{C}$ be the subset of all nets in $\mathcal{D}$ that Cauchy converge in $\mathbb{Q}$. By Lemma \ref{Lem:Nets_R}, it suffices to show that $\langle \rho^{\tilde{\mathcal{Q}}}, \sim^{\tilde{\mathcal{Q}}}, \tilde{+}^{\tilde{\mathcal{Q}}}, \tilde{\times}^{\tilde{\mathcal{Q}}} \rangle \cong \langle \mathcal{C}, \sim, +, \times \rangle$, and that for each $r \in \mathbb{R}$, there is $D \in \mathcal{C}$, such that $r_D = r$. To establish the latter, we define $D_r \in \mathcal{D}$, for each $r \in \mathbb{R}$, as follows:
\[
\textrm{For each state $s$, } D_r(s) =_\df \textrm{the least $q \in s$, such that $|r - q| = \mathrm{min}\{|p - r| : p \in s\}$.}
\]
To see that $r_{D_r} = r$, i.e. that $D_r$ converges to $r$, let $s$ be any state, and suppose $z$ is a positive rational. Let $q$ be a rational number such that $|r - q| < z$. Since $\mathcal{Q}$ is a development model for $\mathbb{Q}$, there is a state $s'$, such that $q \in s'$ and $s' \supseteq s$. Now by definition of $D_r$, we have for all $s'' \supseteq s'$, that $|r - D_r(s'')| \leq z$, as desired. Moreover, since $D_r$ converges to $r$, it is Cauchy convergent in $\mathbb{Q}$. So $D_r \in \mathcal{C}$, as required.

We proceed to establish $\langle \rho^{\tilde{\mathcal{Q}}}, \sim^{\tilde{\mathcal{Q}}}, \tilde{+}^{\tilde{\mathcal{Q}}}, \tilde{\times}^{\tilde{\mathcal{Q}}} \rangle \cong \langle \mathcal{C}, \sim, +, \times \rangle$. The proposed isomorphism, $I$, is defined by $\delta \mapsto D_\delta$, for each $\delta \in \rho^{\tilde{\mathcal{Q}}}$. Let $s$ be a state and let $a, b, c \in \mathrm{Stat}^s$. By Theorem \ref{Thm:Sem_Mirror}:
\begin{align*}
& \tilde{\mathcal{Q}}, s \models (|a - b| < c)^\Diamond  \Longleftrightarrow |a - b| < c \tag{$\dagger$} %\\
%& \tilde{\mathcal{Q}}, s \models (|a - b| < c)^\Diamond \rightarrow \Box \big( (|x - y| < z)^\Diamond \big) \tag{$\ddagger$}
\end{align*}

It follows routinely from ($\dagger$) and the definitions that for any $\xi, \upsilon, \zeta \in \mathrm{Dyn}^{\tilde{\mathcal{Q}}}$:
\begin{align}
& \xi \in \rho^{\tilde{\mathcal{Q}}} \Longleftrightarrow D_\xi \textrm{ Cauchy converges in } \mathbb{Q} \label{Eq:Cauchy_Con} \\
& \xi \sim^{\tilde{\mathcal{Q}}} \upsilon \Longleftrightarrow D_\xi \sim D_\upsilon \label{Eq:Cauchy_Eq} \\
& \tilde{+}^{\tilde{\mathcal{Q}}}(\xi, \upsilon, \zeta) \Longleftrightarrow +(D_\xi, D_\upsilon, D_\zeta) \label{Eq:Add} \\
& \tilde{\times}^{\tilde{\mathcal{Q}}}(\xi, \upsilon, \zeta) \Longleftrightarrow \times(D_\xi, D_\upsilon, D_\zeta)  \label{Eq:Mult}
\end{align}
We establish (\ref{Eq:Cauchy_Con}), which has the most complicated proof of the above. The proof of (\ref{Eq:Cauchy_Eq}) is similar and left the reader. (\ref{Eq:Add}) and (\ref{Eq:Mult}) are immediate from the definitions. Let $\xi \in \mathrm{Dyn}^{\tilde{\mathcal{Q}}}$, and observe the following equivalences:
\begin{align*}
& \tilde{\mathcal{Q}} \models \Box \forall z > 0 \dt \Diamond \Box \exists x \dt \Big( x \leftleftarrows \xi \wedge \Box \big( (|x - \xi| < z)^\Diamond \big) \Big) \\
\Longleftrightarrow & \forall z \in \mathbb{Q}_+ \dt \forall s \in \mathrm{states}_z^{\tilde{\mathcal{Q}}} \dt \exists t \geq^{\tilde{\mathcal{Q}}} s \dt \forall t' \geq^{\tilde{\mathcal{Q}}} t \dt \forall t'' \geq^{\tilde{\mathcal{Q}}} t' \dt \big[ {\tilde{\mathcal{Q}}}, t'' \models \big( |D_\xi(t') - D_\xi(t'')| < z \big)^\Diamond \big] \\
\Longleftrightarrow & \forall z \in \mathbb{Q}_+ \dt \forall s \in \mathrm{states}_z^{\tilde{\mathcal{Q}}} \dt \exists t \geq^{\tilde{\mathcal{Q}}} s \dt \forall t' \geq^{\tilde{\mathcal{Q}}} t \dt \forall t'' \geq^{\tilde{\mathcal{Q}}} t' \dt \big( |D_\xi(t') - D_\xi(t'')| < z \big) \\
\Longleftrightarrow & \forall z \in \mathbb{Q}_+ \dt \exists t \in \mathrm{states}_z^{\tilde{\mathcal{Q}}} \dt \forall t' \geq^{\tilde{\mathcal{Q}}} t \dt \forall t'' \geq^{\tilde{\mathcal{Q}}} t' \dt \big( |D_\xi(t') - D_\xi(t'')| < z \big) \\
\Longleftrightarrow & \forall z \in \mathbb{Q}_+ \dt \exists t \in \mathrm{states}_z^{\tilde{\mathcal{Q}}} \dt \forall t', t'' \geq^{\tilde{\mathcal{Q}}} t \dt \big( |D_\xi(t') - D_\xi(t'')| < z \big)
\end{align*}
The first equivalence follows from the semantics of development models. The second follows from ($\dagger$). The third follows from that $\leq^{\tilde{\mathcal{Q}}}$ is directed. The $\Longleftarrow$ direction of the fourth and last equivalence is immediate. Its $\Longrightarrow$ direction is established as follows: Let $z \in \mathbb{Q}_+$, let $t$ be a state with $z \in t$, and let $t', t'' \geq^{\tilde{\mathcal{Q}}} t$. Since $\geq^{\tilde{\mathcal{Q}}}$ is directed, there is $t''' \geq^{\tilde{\mathcal{Q}}}  t', t''$, such that, by the left-hand-side, $|D_\xi(t') - D_\xi(t''')| < z/2$ and $|D_\xi(t'') - D_\xi(t''')| < z/2$. So $|D_\xi(t') - D_\xi(t'')| < z$, as desired.

By (\ref{Eq:Cauchy_Con}) and $\forall D \in \mathcal{D} \dt (D_{\delta_D} = D)$, we have that $I$ is a bijection from $\rho^{\tilde{\mathcal{Q}}}$ to $\mathcal{C}$. Moreover, by (\ref{Eq:Cauchy_Eq}), (\ref{Eq:Add}), and (\ref{Eq:Mult}), $I$ preserves all the relevant structure. So $I$ is an isomorphism from $\rho^{\tilde{\mathcal{Q}}}$ to $\mathcal{C}$, as desired.
\end{proof}

In terms of the philosophical characterization of development processes in \S \ref{Sec:Char_Dev_Proc}, we can summarize the above account of the reals as follows:

\begin{description}
	\item[States] The set of states can be taken to be any directed partial order, such that for each $r \in \mathbb{R}$, there is a net on this partial order into $\mathbb{Q}$ that Cauchy converges to $r$. Naturally, $\omega$ will do. But when generalizing to completions of other topological spaces, other partial orders may be required.
	\item[Static entities] The static entities are the rational numbers.
	\item[Dynamic entities] The dynamic entities are dynamic rational numbers.
	\item[Type coherence] At any state the dynamic rational manifests as a rational. Thus, the mathematical toolbox for the type of rational numbers applies.
	\item[Teleological development] Any dynamic rational that Cauchy converges may be taken to represent a real number. The teleology of such a dynamic rational can be taken to be a collection of increasingly stringent conditions specifying that the error bound to the real number being represented approaches $0$.
	\item[Non-finality] For any irrational number, the dynamic rational representing it will not satisfy the whole teleology at any state.
\end{description}

\section{Set-theoretic forcing extensions as dynamic $\in$-structures}\label{Sec:Forcing}

In this section we utilize the development of set-theoretic forcing in \cite[Ch. IV.2]{Kun13}, and show that the $\in$-symbol of a forcing extensions can be naturally regarded as a dynamic predicate of a development model. We assume familiarity with basic notions of set theory, defined e.g. in \cite{Kun13}. Assume that $L$ is the language of set theory, with a single relation symbol ${\in} : (\mathrm{Stat}, \mathrm{Stat})$. Let $u$ be a countable transitive model of $\ZF-P$, i.e. $\ZF$ minus the Power set axiom. Let $\mathbb{P} \in u$ be a partially ordered set with ordering $\leq^\mathbb{P}$ and a bottom element $0^\mathbb{P}$. Note that by absoluteness of model-theoretic satisfaction and transitivity of $u$, for any $\phi$ in the language of partial orders, we have $\mathbb{P} \models \phi$ iff $u \models$ ``$\mathbb{P} \models \phi$''. Let $I$ be an arbitrary ideal on $\mathbb{P}$. (N.B.: To fit with the exposition in this paper, the ordering of $\mathbb{P}$ is reversed from the ordering in \citep{Kun13}. Thus,  $I$ corresponds to an arbitrary filter, and $0^\mathbb{P}$ corresponds to a top element, in the development of \cite{Kun13}.)

Following Definitions IV.2.5-6 in \cite{Kun13}, let $u^\mathbb{P}$ be the set of all $\mathbb{P}$-names $\tau$, such that $\tau \in u$. Let $p \in \mathbb{P}$ and $\tau \in u^\mathbb{P}$ be arbitrary. 
%$I_p$ denotes the ideal $\{q \in \mathbb{P} : q \leq p\}$ on $\mathbb{P}$. 
In analogy with Definition IV.2.7 in \cite{Kun13}, we $\in$-recursively define:
\begin{align*}
\mathrm{val}_I(\tau) &=_\df \big\{ \mathrm{val}_I(\sigma) : \exists q \in I \dt (\langle \sigma, q \rangle \in \tau) \big\} \tag{$\dagger$} \\
u[I] &=_\df \big\langle \{ \mathrm{val}_I(\sigma) : \sigma \in u^\mathbb{P} \}, {\in} \restriction u \big\rangle
\end{align*}
%Note that $\mathrm{val}(\tau, I_p) = \big\{ \mathrm{val}(\sigma, I_p) :  \exists q \leq^\mathbb{P} p \dt (\langle \sigma, q \rangle \in \tau) \big\}$.

$u[I]$ is called {\it the forcing extension of $u$ by $I$}. By Lemmata IV.2.10 and IV.2.15 in \cite{Kun13}, $u$ is a substructure $u[I]$, and $u[I]$ is a transitive model of the axioms of Extensionality, Foundation, Pairing, and Union. Moreover, by Theorem IV.2.27 in \cite{Kun13}, if $I$ is a generic ideal, then $u \models \ZF \Longrightarrow u[I] \models \ZF$, and $u \models \ZFC \Longrightarrow u[I] \models \ZFC$.

Assume that $L_\mathrm{DynPred}$ only extends $L$ with the dynamic binary predicate $\dynin$. Now, we construct a natural development model $\mathcal{M}$ of $L_\mathrm{DynPred}^\Diamond$, which has the set of $\mathbb{P}$-names as static individuals, and such that the forcing extension $u[I]$ is isomorphic to a quotient of a definable structure on the static individuals of $\mathcal{M}$. The frame of $\mathcal{M}$ is defined as the restriction of $(\mathbb{P}, \leq^\mathbb{P})$ to $I$, i.e. as follows:
\begin{align*}
\mathrm{states}^\mathcal{M} &=_\df I \\
\leq^\mathcal{M} &=_\df {\leq^\mathbb{P}} \cap I \times I
\end{align*}
Let $s \in I$ be arbitrary. Then $\mathrm{struct}^\mathcal{M}(s)$ is defined as the unique $L_\mathrm{DynPred}$-structure $v$, such that:
\begin{align*}
v \restriction L &= u^\mathbb{P} \\
\mathrm{struct}^\mathcal{M}(s) \models x\dynin y &\Longleftrightarrow x, y \in u^\mathbb{P} \wedge \exists s' \leq^\mathbb{P} s \dt \big( \langle x, s' \rangle \in y \big)
\end{align*}
Note that $\leq^\mathcal{M}$ is directed, since $I$ is an ideal. Therefore, $\mathcal{M}$ is a development model.

Now consider the following $L_\mathrm{Dyn}^\Diamond$-formula expressing the teleological version of $\dynin$:
\begin{align*}
x \tildin y &\leftrightarrow_\df \Diamond \Box \dt x \dynin y
\end{align*}
Note that by the interpretation of $\dynin$ in $\mathcal{M}$, we have for all $\sigma, \tau \in u^\mathbb{P}$: 
\begin{align*}
\mathcal{M} \models \sigma \tildin \tau &\Longleftrightarrow \mathcal{M} \models \Diamond \dt \sigma \dynin \tau \\
&\Longleftrightarrow \exists s \in I \dt \big( \langle \sigma, s \rangle \in \tau \big)
\end{align*}

\begin{thm}
The Mostovsky collapse of $\langle u^\mathbb{P}, \tildin^\mathcal{M} \rangle$ equals the forcing extension $u[I]$.
\end{thm}
\begin{proof}
First, we verify that $\langle u^\mathbb{P}, \tildin^\mathcal{M} \rangle$ is well-founded. Suppose that $\tau_0, \tau_1, \cdots$ is an $\tildin^\mathcal{M}$-descending sequence of length $\omega$. Then there exists an $\omega$-sequence $s_0, s_1, \cdots$ of elements of $I$, such that for all $i < \omega$, $\langle \tau_{i+1}, s_{i+1} \rangle \in \tau_i$. But this contradicts well-foundedness of $\in$.

The Mostovsky collapse, $\mathrm{mos}$, of $\langle u^\mathbb{P}, \tildin^\mathcal{M} \rangle$ is defined by $\tildin^\mathcal{M}$-recursion as $\forall \tau \in u^\mathbb{P} \dt \Big( \mathrm{mos}(\tau) = \{ \mathrm{mos}(\sigma) : \sigma \tildin^\mathcal{M} \tau \} \Big)$. So by the remark preceding this theorem, we have 
\begin{align*}
\forall \tau \in u^\mathbb{P} \dt \Big( \mathrm{mos}(\tau) = \Big\{ \mathrm{mos}(\sigma) : \exists s \in I \dt \big( \langle \sigma, s \rangle \in \tau \big) \Big\} \Big).
\end{align*}
Thus, the definition of $\mathrm{mos}$ is the same as the definition ($\dagger$) of $\mathrm{val}_I$. So the Mostovsky collapse of $\langle u^\mathbb{P}, \tildin^\mathcal{M} \rangle$ equals the forcing extension $u[I]$, as desired. (In particular, $\mathrm{mos} : u^\mathbb{P} \rightarrow u[I]$ is onto, and for any $\sigma, \tau \in u^\mathbb{P}$, we have $\sigma \tildin^\mathcal{M} \tau \rightarrow \mathrm{mos}(\sigma) \in \mathrm{mos}(\tau)$.)
\end{proof}

We can regard the Mostovsky collapse of $\langle u^\mathbb{P}, \tildin^\mathcal{M} \rangle$ as what you get by quotienting $\langle u^\mathbb{P}, \tildin^\mathcal{M} \rangle$ by the most fine-grained equivalence relation that renders the quotient a model of extensionality. 

\begin{que*}
Does $\langle u^\mathbb{P}, \tildin^\mathcal{M} \rangle$ define an equivalence relation $x \sim y$ on $u^\mathbb{P}$, such that for all $\tau, \tau' \in u^\mathbb{P}$, $\big[ \langle u^\mathbb{P}, \tildin^\mathcal{M} \rangle \models \tau \sim \tau' \big] \Longleftrightarrow \mathrm{mos}(\tau) = \mathrm{mos}(\tau')$?
\end{que*}

The development model $\langle u^\mathbb{P}, \tildin^\mathcal{M} \rangle$ also invites us to view the elements of the forcing extension as dynamic elements. Indeed, for any $\tau \in u^\mathbb{P}$, consider the functional dynamic individual $D_\tau$, defined as follows, for each $s \in I$: 
\begin{align*}
D_\tau(s) = \big\{ \langle \sigma, 0^\mathbb{P} \rangle : \exists q \leq^\mathbb{P} s \dt (\langle \sigma, q \rangle \in \tau) \big\} \in u^\mathbb{P}
\end{align*}
(In the above set-abstract, think of $\langle \sigma, 0^\mathbb{P} \rangle$ as encoding $\sigma$. The point of using $\langle \sigma, 0^\mathbb{P} \rangle$ is that it ensures that $D_\tau(s) \in u^\mathbb{P}$, for all $s$.) Note that for any $s \in I$, and any $\sigma, \tau \in u^\mathbb{P}$, we have $\langle \sigma, 0 \rangle \in D_\tau(s) \Longleftrightarrow \mathcal{M}, s \models \sigma \dynin \tau$. Moreover, for any $\sigma, \tau \in u^\mathbb{P}$, we have $\exists s \in I \dt \big( \langle \sigma, 0 \rangle \in D_\tau(s) \big) \Longleftrightarrow \mathcal{M}, s \models \sigma \tildin \tau$. So for any $\sigma, \tau \in u^\mathbb{P}$, we have intuitively that $D_\tau$ teleologically has $\sigma$ as element iff $\mathrm{val}_I(\tau)$ has $\mathrm{val}_I(\sigma)$ as element (in the forcing extension $u[I]$).

\begin{description}
	\item[States] The set of states can be taken to be any ideal, $I$, on any partial order, $\mathbb{P}$.
	\item[Static entities] The static entities are the $\mathbb{P}$-names.
	\item[Dynamic entities] The dynamic entity is the dynamic membership-relation.
	\item[Type coherence] At any state, the dynamic membership-relation is definable in the ground model, and can thus be managed with the mathematical resources of the ground model.
	\item[Teleological development] The teleology of the dynamic membership-relation consists of the atomic facts of the membership-relation of the forcing extension.
	\item[Non-finality] In general, the ideal is not definable in the ground model. In such cases, the whole teleology will not be met at any particular state. 
\end{description}

\section{Realizers of types and non-standard numbers as dynamic }\label{Sec:Saturation}

Assume (for ease of presentation) that $L$ is countable,\footnote{This section generalizes immediately to uncountable languages.} and consider an $L$-model $u$.
% $\mathbb{N}$ of arithmetic, in the language $L$, here taken to be the purely relational language $L_\PA$ introduced in \S \ref{Sec:Background}. 
Let $p(\vec{x})$ be an enumerated set of $L$-formulas, with free variables $\vec{x}$. Then we write $p_0(\vec{x}), p_1(\vec{x}), p_2(\vec{x}), \cdots$ for the enumeration. $p(\vec{x})$ is a {\it  type over $u$} if for any finite $\Gamma \subseteq p(\vec{x})$, there is $\vec{a} \in \mathbb{N}$, such that for each $\phi(\vec{x}) \in \Gamma$, we have $u \models \phi(\vec{a})$. Assume that $p(\vec{x})$ is a type. $p(\vec{x})$ is {\it  complete} if for any $\phi(\vec{x}) \in L$, $\phi \in p$ or $\neg\phi \in p$.  $p(\vec{x})$ is {\it  realized in} $u$ if there is $\vec{a} \in u$, such that $u \models p(\vec{a})$.

Let $\mathbb{A}$ be a set of finite subsets of $\omega$, such that $\langle \mathbb{A}, \subseteq \rangle$ is directed and $\bigcup \mathbb{A} = \omega$ (e.g. we could have $\mathbb{A}$ be the set of finite initial segments of $\omega$), and let $\mathcal{M}$ be the simple development model for $u$ over $(\mathbb{A}, \subseteq)$. Let $\mathcal{P}$ be some set of types over $u$. For each $p(\vec{x}) \in \mathcal{P}$, there is a functional dynamic tuple $D_p : \omega \rightarrow u^{\lng(\vec{x})}$, such that for all $F \in \mathbb{A}$, we have $u \models \bigwedge_{i \in F} p_i(D_p(i))$. Let $\mathcal{D} = \{ D_p : p \in \mathcal{P}\}$. The least development model extending $\mathcal{M}$ by $\mathcal{D}$, denoted $\tilde{\mathcal{M}}(\mathcal{D})$,\footnote{Recall the definition from the end of \S \ref{Sec:Prel}} formalizes a development process encapsulating all such dynamic tuples. For each $p \in \mathcal{P}$, the dynamic tuple in $\mathrm{Dyn}^{\tilde{\mathcal{M}}(\mathcal{D})}$ corresponding to $D_p$ is denoted $\delta_p$. Note that for each $p \in \mathcal{P}$, we have $\tilde{\mathcal{M}}(\mathcal{D}) \models^\mathrm{tele} p(\delta_p)$.\footnote{This is also stated and elaborated in the statement and proof of Theorem \ref{Thm:Dynamic_Non-standard}.} This shows that for each type $p \in \mathcal{P}$, the dynamic tuple $\delta_p$ realizes $p$ in a teleological sense.

A technique from model theory for realizing types over $u$ is to form the ultrapower $u^\omega / \mathcal{U}$ of $u$, for some non-principal ultrafilter $\mathcal{U}$ on $\omega$. Its elements are the equivalence classes, $[f]_\mathcal{U}$, of functions $f : \omega \rightarrow u$, given by an equivalence relation, $\sim_\mathcal{U}$ on $u^\omega$, defined by $f \sim_\mathcal{U} g \Longleftrightarrow_\df \{x \in \omega : f(x) = g(x)\} \in \mathcal{U}$. The last fact we need about ultrapowers is the following classical theorem:
\begin{thm}[\L{}o\'s]
Let $u, \mathcal{U}$ be as above. Let $\phi(x_1, \cdots, x_n) \in L$, and let $f_1, \cdots, f_n$ be functions from $\omega$ to $u$.
\begin{align*}
u^\omega / \mathcal{U} \models \phi([f_1]_\mathcal{U}, \cdots, [f_n]_\mathcal{U}) \Longleftrightarrow \big\{ s < \omega : u \models \phi(f_1(s), \cdots, f_n(s)) \big\} \in \mathcal{U}
\end{align*} 
\end{thm}

Moreover, assume that $\mathbb{A}$ is the set of finite initial segments of $\omega$. Since $\langle \mathbb{A}, \subseteq \rangle$ is isomorphic to $\langle \omega, \leq^\omega \rangle$, we may identify their elements by this (unique) isomorphism. Let $p(\vec{x})$ be a type over $u$. Note that under the above identification, the functional dynamic individual of $p$ is a function $D_p : \omega \rightarrow u$, i.e. it is an element of $u^\omega$. Thus, the ultrapower $u^\omega / \mathcal{U}$ is a quotient of the structure on the set of all functional dynamic individuals of $\mathcal{M}$, obtained by interpreting the primitive arithmetic relations statewise.

By a natural and minor extension of notation, we define $[D_p]_\mathcal{U} =_\df \big\langle [\pi_1 \circ D_p]_\mathcal{U}, \cdots, [\pi_{\lng(\vec{x})} \circ D_p]_\mathcal{U} \big\rangle$. The following theorem establishes a close correspondence between the dynamic tuple $\delta_p$ and the possibly non-standard number $[D_p]_\mathcal{U}$:

\begin{thm}\label{Thm:Dynamic_Non-standard}
Let $\mathcal{P}$, $\tilde{\mathcal{M}}(\mathcal{D})$, and $u^\omega / \mathcal{U}$ be as above. Let $p(\vec{x}) \in \mathcal{P}$ and $\phi(\vec{x}) \in L$.
\begin{align*}
\phi \in p \Longrightarrow \tilde{\mathcal{M}}(\mathcal{D}) \models^\mathrm{tele} \phi(\delta_p) \Longrightarrow u^\omega / \mathcal{U} \models \phi([D_p]_\mathcal{U}) 
\end{align*}
Moreover, if $p$ is complete, then:
\begin{align*}
\phi \in p \Longleftrightarrow \tilde{\mathcal{M}}(\mathcal{D}) \models^\mathrm{tele} \phi(\delta_p) \Longleftrightarrow u^\omega / \mathcal{U} \models \phi([D_p]_\mathcal{U}) 
\end{align*}
\end{thm}
\begin{proof}
The first $\Longrightarrow$: Suppose that $\phi \in p$. Then $\phi$ is $p_s$, for some $s \in \omega$. So by construction of $\delta_p$, we have for all $t \geq s$ that $\tilde{\mathcal{M}}(\mathcal{D}), t \models \phi(\delta_p)$. Therefore, $\tilde{\mathcal{M}}(\mathcal{D}) \models \Diamond \Box \phi(\delta_p)$. 

The last $\Longrightarrow$: Suppose that $\tilde{\mathcal{M}}(\mathcal{D}) \models \Diamond \Box \phi(\delta_p)$. Then there is $s \in \omega$, such that for all $t \geq s$, we have $\tilde{\mathcal{M}}(\mathcal{D}), t \models \phi(\delta_p)$. Since $\mathcal{U}$ is non-principal, it has $\{t \in \omega : t \geq s\}$ as an element. So $u^\omega / \mathcal{U} \models \phi([D_p]_\mathcal{U})$.

The last two $\Longleftrightarrow$: Assume now that $p$ is complete, and that $u^\omega / \mathcal{U} \models \phi([D_p]_\mathcal{U})$. If $\phi \not\in p$, then $\neg\phi \in p$. So by the above, $u^\omega / \mathcal{U} \models \neg\phi([D_p]_\mathcal{U})$, a contradiction. Hence, $\phi \in p$.
\end{proof}

We shall now apply the above to arithmetic, i.e. to the case that $L$ is the relational language of arithmetic, and $u$ is the standard model $\mathbb{N}$, considered as an $L$-structure. Note that the non-standard model $\mathbb{N}^\omega / \mathcal{U}$ is a quotient of the statewise interpretation of the arithmetic relations on the set of all functional dynamic individuals. Assume for simplicity that $\mathcal{P}$ is the set of all complete types over $u$.\footnote{The case that $\mathcal{P}$ is the set of all recursive types is also highly significant, but out of scope for the present paper.} Then, by Theorem \ref{Thm:Dynamic_Non-standard}, we obtain for each $\phi \in L$ and each complete type $p$ over $u$:
\begin{align*}
\phi \in p \Longleftrightarrow \tilde{\mathcal{M}}(\mathcal{D}) \models^\mathrm{tele} \phi(\delta_p) \Longleftrightarrow \mathbb{N}^\omega / \mathcal{U} \models \phi([D_p]_\mathcal{U}) 
\end{align*}
Thus, the dynamic numbers are representatives of non-standard (and standard) numbers. Moreover, the properties of non-standard numbers in an ultrapower of $\mathbb{N}$ are mirrored by the teleological properties of the corresponding dynamic numbers in a natural development model for $\mathbb{N}$. We proceed to connect this back to the philosophical development of \S \ref{Sec:Char_Dev_Proc}:

\begin{description}
\item[States] The set of states can simply be taken to be $\omega$ with its usual order (but it naturally generalizes to take the finite subsets of $\omega$ as the states). The structure of each state $s$ is an expansion of $\mathbb{N}$, accounting for the manifestation of each dynamic entity at $s$.
\item[Static entities] The static entities are the standard natural numbers.
\item[Dynamic entities] The dynamic entities are dynamic numbers and tuples of numbers: For each complete first-order description,  $p(\vec{x}) \in \mathcal{P}$, of a non-standard number, there is a dynamic tuple $\delta_p$, which manifests at any state $s$ to satisfy $\bigwedge_{j \leq s} p_j(\delta_p)$.
\item[Teleological development] The teleology of each dynamic tuple $\delta_p$ is the collection of conditions formalized by the type $p(\vec{x})$. For any state $s$, any teleological formulation $p(\vec{x})$, and any $p_i(\vec{x}) \in p$, there is the potential to transition to the state $\max\{s, i\}$, at which $\delta_p$ manifests to satisfy $p_i(\vec{x})$ at all accessible states.
\item[Non-finality] For $p(x) = \{x = 0\}$, obviously $\delta_p$ satisfies the whole teleology at every state. But for $p(x) = \{n < x : n \in \mathbb{N}\}$, $\delta_p$ does not satisfy the whole teleology at any state, and $\delta_p$ robustly represents a non-standard number.
\end{description}

\section{A development process of reflection in set theory}\label{Sec:Reflection}

The iterated application of the Reflection Theorem of $\ZF$, as applied in the cumulative hierarchy of $V_\alpha$s, can be construed as a development process with the following characteristics:

\begin{description}
\item[Language] $L$ is the language of set theory.
\item[States] Each state $s$ is an ordered pair $\langle \Gamma_s, \vec{a}_s \rangle$ consisting of a finite set $\Gamma_s$ of $L$-formulas paired with a finite tuple $\vec{a}_s$ of sets. The structures of all states is constant; it is an underlying model $u$ of $\ZF$.
\item[Static entities] The sets.
\item[Dynamic entities] A dynamic structure, $\nu$, which manifests in each state $s$ as the least $V_\alpha$, such that for each $\phi(\vec{x}, \vec{y}) \in \Gamma$, where $\lng(\vec{y}) = \lng(\vec{a})$, we have that $u \models \forall \vec{x} \dt \big( \phi(\vec{x}, \vec{a}) \leftrightarrow \phi(\vec{x}, \vec{a})^\nu \big)$, where $\phi(\vec{x}, \vec{a})^\nu$ is the formula obtained by restricting all quantifiers in $\phi$ to the manifestation of $\nu$. The existence of such a $V_\alpha$ follows from the reflection theorem of $\ZF$.
\item[Teleological development] For any finite $\Gamma \subseteq L$ and any tuple $\vec{a}$ of sets, $\forall \vec{x} \dt \big( \phi(\vec{x}, \vec{a}) \leftrightarrow \phi(\vec{x}, \vec{a})^\nu \big)$ is satisfied teleologically.
\item[Non-finality] For each state $s$, there is a tuple $\vec{a}$ not in the manifestation of $\nu$, so the condition  $\forall \vec{x} \dt \big( \phi(\vec{x}, \vec{a}) \leftrightarrow \phi(\vec{x}, \vec{a})^\nu \big)$ is not met in that state.
\end{description}

Given a model $u \models \ZF$, this development process is captured by a development model $\mathcal{M}$ constructed as outlined above. We obtain for each $\phi(\vec{x}, \vec{y}) \in L$ and each $\vec{a} \in u$, where $\lng(\vec{y}) = \lng(\vec{a})$:
\begin{align*}
\mathcal{M} \models^\mathrm{tele} \forall \vec{x} \dt \big( \phi(\vec{x}, \vec{a}) \leftrightarrow \phi(\vec{x}, \vec{a})^\nu \big)
\end{align*}

\section{A development process of truth revision}\label{Sec:Revision_semantics}

The liar paradox, attributed to Epimenides of Crete ($6^\textrm{th}$ BC), can be formalized (along the lines of Tarski's \cite{Tar36}) as that the following axiom schema over $\PA$ is inconsistent: For each sentence $\phi$ of the functional language of $\PA$ augmented with a truth predicate $T$: 
\begin{align*}
T(\gquote{\phi}) \leftrightarrow \phi \tag{T}
\end{align*}

We shall follow the formalism of \S \ref{Sec:Prel}, except that $L$ is allowed to be functional.\footnote{The reason that we work with a functional language in this section is that it is convenient to have terms referring to syntax in the language when reasoning about truth. As explained below, this deviation from the assumption made in \S \ref{Sec:Prel} (that all languages be relational) does not pose any difficulty in the present context.} Assume that $L$ is the usual language of $\PA$ with the constant and function symbols $\underline{0}, S, +, \times$, but formalized with a single sort $\mathrm{Stat}$. Now assume that $\PA$ is formalized in the language $L_T$ of arithmetic obtained by augmenting $L$ with the truth predicate $T(x)$, intended to formalize ``$x$ is true''. Note that the induction schema of $\PA$ does not have instances for formulas with the truth predicate; let $\PA(L_T)$ be $\PA$ with the induction schema extended to all formulas of $L_T$.\footnote{The inconsistency of (T) does not require this extended induction schema. However, it is philosophically natural, so it is desirable that a theory of truth be consistent with it.} For every formula $\phi$ in this language, there is a so called G\"odel code $\gquote{\phi}$, which is a term in the language (intended to refer to the formula $\phi$). So if $\phi$ is a sentence, then $T(\gquote{\phi})$ formalizes ``$\phi$ is true''.

The inconsistency of (T) over $\PA$ stems from a single sentence $\lambda$, which formalizes a liar sentence expressing ``This sentence is not true.'' in the sense that $\PA \vdash \neg T(\gquote{\lambda}) \leftrightarrow \lambda$. The existence of such a $\lambda$ follows from G\"odel's fixed point lemma. Today there is a rich landscape of consistent theories of truth with axiomatizations side-stepping the full schema (T) above, see \cite{Hal14}. One of the most prominent of these is the Friedman--Sheard system ($\FS$) from  
\cite{FS87}. 

$\FS$ has several philosophically significant and desirable features: Firstly, in \cite{Hal94} $\FS$ was given a satisfying axiomatization reflecting the Tarskian inductive definition of truth from \cite{Tar36}. Secondly, $\FS$ proves that the truths adhere to classical logic. Thirdly, $\FS$ treats external and internal truth symmetrically in the sense that it proves $\phi$ iff it proves ``$\phi$ is true''. Fourthly, $\FS$ proves a theorem schema, which intuitively is the restriction of (T) to its typable instances, and which thereby provides a deep explanation of the liar paradox; see \cite[Corollary 14.24]{Hal14}. The main philosophical drawback of $\FS$ as a theory of truth is that it is $\omega$-inconsistent, and hence only has non-standard models. This follows from McGee's paradox \cite{McG85}. 

As a vindication of $\FS$, this section shows that $\FS$ is teleologically satisfied by a development model that has a standard model of arithmetic at each state, and treats $T$ as a dynamic predicate whose extension changes as the process unfolds. Thus, $\FS$ can be philosophically endorsed without committing to a non-standard model of arithmetic, as long as one is willing to regard truth as dynamic in the sense of this paper. This result relies on that $\FS$ is locally validated (in a precise sense explained below) by the revision semantics of truth, as shown in \cite{Hal94}. The revision semantics was independently introduced by Gupta and Herzberger in \cite{Gup82} and \cite{Her82a, Her82b}, respectively. %This section shows that this revision semantics may naturally be viewed as a development process in the sense explained in \S \ref{Sec:Char_Dev_Proc} and formalized in \S \ref{Sec:Prel}. 

Since we use a functional language in this section, the definition of development model will be slightly tweaked. The problem with having function symbols in a modal setting is that a particular world may not be closed under the functions. However, in this section all the structures of the states are expansions of $\mathbb{N}$, so they are all closed under the functions considered. Therefore, this tweaked definition does not pose any difficulties.

In order to discuss truth as a formal predicate, some notation concerning syntax needs to be introduced. For any arity $n \in \mathbb{N}$:
\begin{align*}
\Var &=_\df \textnormal{the set of variables.} \\
\Term_{L_T}^n &=_\df \textnormal{the set of $L_T$-terms with precisely $n$ free variables.} \\
L_T^n &=_\df \textnormal{the set of $L_T$-formulas with precisely $n$ free variables.} \\
L_T^n[v_1, \cdots, v_n] &=_\df \textnormal{the set of $L_T$-formulas with precisely $v_1, \cdots, v_n$ as free variables.} 
%\Atom_L^n &=_\df \textnormal{the set of atomic $L$-formulas with precisely $n$ free variables.} 
\end{align*}
When the superscript $n$ above is dropped, all arities are included, e.g. $\Term_{L_T}$ is the set of all $L_T$-terms. The {\it  substitution function}, denoted $\sbt(x, y, z)$, maps $\langle \phi, t, v \rangle$ to the formula obtained from $\phi$ by substituting every instance of the variable $v$ by the term $t$. For any $n \in \mathbb{N}$, $\underline{n}$ denotes the {\it  numeral} $S^n(\underline{0}) \in \Term_{L_T}^0$. Let $t \in \Term_{L_T}$ and $\phi \in {L_T}$. The {\it  G\"odel quotes} $\gq(t)$ (or $\gquote{t}$) and $\gq(\phi)$ (or $\gquote{\phi}$) denote numerals which represent $t$ and $\phi$, respectively, in arithmetic. Closely related to these are the {\it  G\"odel codes} $\mathrm{gc}(t), \mathrm{gc}(\phi)$ of $t$ and $\phi$, respectively, which are the numbers in the standard model interpreting $\gq(t), \gq(\phi)$, respectively. Moreover, $\PA$ is assumed to be equipped with various function and relation symbols, representing various syntactically relevant functions and relations.\footnote{Formally, these symbols could alternatively be implemented by means of a single relation symbol.} These representations are denoted with a dot under the name of the relation/function in question. For example, $\udot L_T^0$ denotes the representation of $L_T^0$, and $\udot \wedge$ represents a function from $L_T \times L_T$ to $L_T$ mapping two formulas to their conjunction. The evaluation function returning the numeric value of any closed $L_T$-term is denoted $\ev : \Term_{L_T}^0 \rightarrow \mathbb{N}$. Moreover, for any recursive system $S$, there is a formula $\Pr_S$ in $L$ that {\it  adequately represents provability in $S$}, in the precise sense that it satisfies the well-known Hilbert--Bernays provability conditions. With the above notation in place, $\FS$ can be axiomatized as in Figure \ref{Fig:FS}.

\begin{figure}

\caption{Axioms and rules of $\FS $}
\label{Fig:FS}

\begin{center}
 
\vspace{12pt}
{\bf The system $\FS $}

\vspace{12pt}
$\FS$ extends $\PA$ with the following axioms and rules:

\vspace{18pt}
{\bf Axioms of $\FS $}
\end{center}
\begin{align*}
\forall t_1, t_2 \in \udot \Term_{\udot{L}_T}^0 & \dt \big( T \big( \gquote{t_0 = t_1} \big)  \leftrightarrow \udot\ev(t_0) = \udot\ev(t_1) \Big) \tag{$\UCT^\Atom$}  \\
 \forall \phi \in \udot L_T^0 & \dt \big( T(\udot \neg \phi)   \leftrightarrow \neg T(\phi) \big) \tag{$\UCT^\neg$} \\
 \forall \phi, \psi \in \udot L_T^0 & \dt \big( T(\phi \udot\rightarrow \psi)   \leftrightarrow (T(\phi) \rightarrow T(\psi)) \big) \tag{$\UCT^\rightarrow$} \\
 \forall \phi, \psi \in \udot L_T^0 & \dt \big( T(\phi \udot\wedge \psi)   \leftrightarrow (T(\phi) \wedge T(\psi)) \big) \tag{$\UCT^\wedge$} \\
 \forall \phi, \psi \in \udot L_T^0 & \dt \big( T(\phi \udot\vee \psi)   \leftrightarrow (T(\phi) \vee T(\psi)) \big) \tag{$\UCT^\vee$} \\
 \forall v \in \udot \Var \dt \forall \phi \in \udot L_T^1[v] & \dt \big( T(\udot \forall v \phi)   \leftrightarrow \forall t \in \udot\Term_{\udot {L_T}}^0 \big( T(\udot\sbt(\phi, t, v)) \big) \Big) \tag{$\UCT^\forall$} \\
 \forall v \in \udot \Var \dt \forall \phi \in \udot L_T^1[v] & \dt \big( T(\udot \exists v \phi)   \leftrightarrow \exists t \in \udot\Term_{\udot {L_T}}^0 \big( T(\udot\sbt(\phi, t, v)) \big) \Big)  \tag{$\UCT^\exists$} 
\end{align*}

\vspace{6pt}
\begin{center} 
{\bf Rules of $\FS $}
\end{center}
\begin{align*}
\vdash \phi \phantom{~~} &\Longrightarrow \phantom{~~} \vdash T(\gquote{\phi}) \textrm{, for each $\phi \in L_T^0$.}  \tag{$\NEC$} \\
\vdash \phi \phantom{~~} &\Longleftarrow \phantom{~~} \vdash T(\gquote{\phi}) \textrm{, for each $\phi \in L_T^0$.}  \tag{$\CONEC$}
\end{align*}

\end{figure}

The revision semantics starts with a ground model $\langle \mathbb{N}, P_0 \rangle$ expanding the standard model of arithmetic with an arbitrary interpretation $P_0 \subseteq \{ \mathrm{gc}(\sigma) : \sigma \in L_T^0 \}$ of the truth predicate (e.g. we may take $P_0$ to be $\varnothing$), and iteratively revises the interpretation of the truth predicate to $P_1, P_2, \cdots$, by the recursion 
\[ P_{n+1} =_\df \big\{ \mathrm{gc}(\sigma) : \langle \mathbb{N}, P_n \rangle \models \sigma \big\}. \]
The steps are obtained by the map $\Gamma$, defined by $P \mapsto \big\{ \mathrm{gc}(\sigma) : \langle \mathbb{N}, P \rangle \models \sigma \big\}$, for all $P \subseteq \mathbb{N}$.

Let $\FS_0$ be $\PA(L_T)$. Let $\FS_1$ be $\FS$ minus the rules $\NEC$ and $\CONEC$, but plus the axiom $\forall \phi \in \udot L_T^0 \big( \Pr_{\PA(L_T)}(\phi) \rightarrow T(\phi) \big)$. For each $2 \leq n \in \mathbb{N}$, let $\FS_n$ be the Hilbert style deductive subsystem of $\FS$, obtained by the restriction that at most $n-1$ applications of $\NEC$ and at most $n-1$ applications of $\CONEC$ are allowed in a proof.

We proceed to explicate the revision semantics as a development model, with a tweak of the definition of development model in \S \ref{Sec:Prel} allowing that the underlying language $L$ has function symbols.\footnote{As explained above, in the present context this deviation is insignificant and safe.} Assume that $L_\mathrm{DynPred}$ is $L$ augmented as in \S \ref{Sec:Prel} with a dynamic predicate $T$. Now we construct a development model $\mathcal{M}$ of $L_\mathrm{DynPred}^\Diamond$ as follows: The frame $\langle \mathrm{states}^\mathcal{M}, \leq^\mathcal{M} \rangle$ is defined to be $\langle \omega, \leq^\omega \rangle$. For each $s \in \mathrm{states}^\mathcal{M}$, $\mathrm{Struct}^\mathcal{M}(s)$ is defined as the unique $L_\mathrm{DynPred}$-model, such that:
\begin{align*}
\mathrm{Struct}^\mathcal{M}(s) \restriction L &= \mathbb{N} \\
a \in T^\mathcal{M} &\Longleftrightarrow a \in P_s
\end{align*}

We expand on the properties of $\mathcal{M}$ in terms of the terminology of \S \ref{Sec:Char_Dev_Proc}:

\begin{description}
\item[States] The states are the natural numbers $\omega$ with their usual ordering as the accessibility relation.
\item[Static entities] The static entities are the natural numbers and ambient standard model $\mathbb{N}$ of arithmetic, which is the underlying structure of all the states.
\item[Dynamic entities] The dynamic predicate $T$ is interpreted in each state $s$ as the class $P_s$ of G\"odel codes of true sentences in that state.
\item[Transformation operation] For any state $s$, the structure of the state $s+1$ is obtained by transforming the structure of $s$, by means of the revision operator: $P_{s+1} = \Gamma(P_s)$, where $P_s, P_{s+1}$ are the interpretations of $T$ at $s, s+1$, respectively.
\item[Teleological development] The teleology consists of the theorems of $\FS$, which $T$ develops in order for the states to satisfy. The following result (proved in \cite[Theorem 14.15]{Hal14}) is key:
\begin{thm}\label{Thm:Rev_Tel_Dev}
For all subsets $Q \subseteq \{ \mathrm{gc}(\sigma) : \sigma \in L_T^0 \}$:
\[ \langle \mathbb{N}, Q \rangle \models \FS_n \Longleftrightarrow \exists Q_0 \subseteq \mathbb{N} \big( Q = \Gamma^n(Q_0) \big) \]
\end{thm}
The relevant consequence for $\mathcal{M}$ is the following:
\begin{cor}
$\mathcal{M} \models^\mathrm{tele} \FS$
\end{cor}
\item[Non-finality] It follows from McGee's theorem \cite{McG85}, that $\FS$ is $\omega$-inconsistent. Hence, it does not have a standard model, and is therefore not fully satisfied at any state of $\mathcal{M}$.
\end{description}

The main philosophical gain of casting the revision semantics as the development model $\mathcal{M}$ is that $\mathcal{M}$ has a standard model of arithmetic at every state, and nonetheless satisfies $\FS$ in the teleological sense. This allows us to philosophically endorse $\FS$ as an adequate theory of truth without committing to non-standard models of arithmetic. It prompts us, on the other hand, to regard the class of truths as a dynamic class which develops teleologically.

Finally, it is worth pointing out how $\mathcal{M}$ handles the liar sentence $\lambda$. Note that for any conventional first-order model $v$ of $\FS$, we have $v \models \lambda$ or $v \models \neg \lambda$, and both options have models. This is somewhat dissatisfying since it introduces an arbitrary choice that we might not want to make. In the model $\mathcal{M}$, the arbitrariness is removed. Instead we obtain $\mathcal{M} \models \Box (\Diamond \lambda \wedge \Diamond \neg \lambda)$, illuminating the dynamics of the liar.\footnote{To avoid arbitrariness with regard to a truth-teller $\tau$, such that $\PA \vdash T(\tau) \leftrightarrow \tau$, we need to pass to a slightly more complex model $\mathcal{M}'$. $\mathcal{M}'$ is like $\mathcal{M}$ except that it has many isolated copies of $\langle \omega, \leq^\omega \rangle$ as its states, such that the bottom states are mapped surjectively by $\mathrm{Struct}^{\mathcal{M}'}$ to all expansions of $\mathbb{N}$ interpreting truth (i.e. we allow arbitrary starting points $P_0$ in the revision semantics). We then avoid $\mathcal{M}' \models \Box \tau$ and $\mathcal{M}' \models \Box \neg \tau$, while obtaining $\mathcal{M}' \models \Box \tau \vee \Box \neg \tau$, again illuminating the dynamics.} This is shown by a routine exercise establishing that if $\lambda \in Q$, for some $Q \subseteq \{ \mathrm{gc}(\sigma) : \sigma \in L_T^0 \}$, then $\neg \lambda \in \Gamma(Q)$. In other words, the revision operator switches between taking $\lambda$ to be true and taking $\neg \lambda$ to be true.

\section{Conclusion}\label{Sec:Conc}

The wide variety of mathematical constructions accounted for as development processes in this paper support that the notion plays a unifying role in the philosophy of mathematics, being both general, natural and expedient. Let us summarize how this plays out in the applications covered:
\begin{itemize}
	\item In the case of the reals we have seen that a topological net into the rational line is the same thing as a (functional) dynamic rational number. This naturally generalizes to other topological spaces, and in light of \cite[p. 187]{Mun00} we are therefore justified in conjecturing that the completion of any Hausdorff topological space may be accounted for as a development process. This shows that the continuum can be understood without appeal to idealized limits or sets of rational numbers, but instead as the development of dynamic rational approximations.
	\item In the case of forcing, the forcing conditions are organized as a directed partial order, and we can consider these as states of a development process. This teleologically accounts for the $\in$-relation of the forcing extension as a dynamic relation. Thus, the seemingly mysterious act of adjoining new sets can be recast as a natural process of development, bringing clarity to the ontology of forcing.
	\item The realization of types in model theory is a notion that is designed to capture a very general and ubiquitous construction in mathematics. So the fact that these can be accounted for as dynamic objects that develop teleologically shows the wide relevance of this procedural view. In particular, this accounts for many non-standard natural numbers as dynamic natural numbers. A major advantage of this approach is that a development model with such dynamic natural numbers has the standard model $\mathbb{N}$ at each state. Thus, it accounts for non-standard numbers by means of the standard structure. This provides a philosophically attractive reconciliation: we can explain the non-standard numbers as dynamic numbers, without inflating the underlying ontology beyond the familiar standard naturals.
	\item The Reflection Theorem of $\ZF$ has a central role in the foundations of mathematics, formalizing the vast height of the iterative conception of sets. We have seen that the reflecting sub-universes can be accounted for as a single dynamic sub-universe that teleologically develops to fully reflect the entire universe. This suggests that the reflection phenomenon is not a brute fact about static set-theoretic height, but the natural outcome of a teleological process that unfolds a dynamic universe of sets.
	\item The revision semantics provides an attractive approach to the liar paradox, as it yields a theory of truth ($\FS$) that is both classical (the truths are closed under classical logic) and untyped (there is only one truth predicate and it applies to itself). However, this theory may be considered unnatural as it is $\omega$-inconsistent, and therefore only has non-standard models. But we have seen that the revision semantics can be modeled by a development model, in which the truth predicate is accounted for as a dynamic predicate. As this account is entirely based on the standard model $\mathbb{N}$, it is much more natural, and thus vindicates the revision semantics. This shows how the perspective of development-processes restores philosophical plausibility to a major theory of truth by grounding it in the standard structure.
\end{itemize}

What these diverse applications ultimately show is the unifying strength of the view proposed: mathematical constructions can be understood as potentialist processes in which dynamic objects develop state by state to satisfy a teleology, while retaining a coherent identity and type. It not only systematizes familiar constructions but also opens a promising avenue for understanding the ontology of infinitary mathematics as radically dynamic.

\section*{Funding}

This research was supported by Vetenskapsr{\aa}det (The Swedish Research Council), grant ID 2020-00613.

\end{document}